\documentclass[reqno]{amsart}
\usepackage[english]{babel}
\usepackage[T1]{fontenc}
\usepackage[utf8]{inputenc}
\usepackage{amsmath}
\usepackage{amsfonts}
\usepackage{array}
\usepackage{eurosym}
\usepackage{amssymb}
\usepackage{mathrsfs}
\usepackage{amsthm}
\usepackage{xfrac}
\usepackage{lmodern}
\usepackage{fix-cm}
\usepackage{listings}
\usepackage{fancyvrb}
\usepackage{enumerate}
\usepackage{xcolor}
\usepackage{pgf}
\usepackage{tikz}
\tikzstyle{every picture}=[line width=.7pt,minimum size=3pt,every label/.append style={font=\normalsize},label distance=2pt]
\tikzstyle{every node}=[font=\normalsize,circle,draw=black,fill=black,inner sep=0pt,minimum width=1.3pt]
\usetikzlibrary{arrows}
\usepackage{mathrsfs}
\usepackage{float}
\usepackage{geometry}
\usepackage{caption}
\usepackage{subcaption}
\usepackage[labelformat=simple]{subcaption}

\usepackage{color}
\usepackage[bookmarks=true]{hyperref}

\usepackage{wasysym}
\usepackage{enumitem}
\usepackage{longtable}

\makeatletter
\newcolumntype{"}{@{\hskip\tabcolsep\vrule width 1pt\hskip\tabcolsep}}
\makeatother

\makeatletter
\newtheorem*{rep@theorem}{\rep@title}
\newcommand{\newreptheorem}[2]{%
\newenvironment{rep#1}[1]{%
 \def\rep@title{#2 \ref{##1}}%
 \begin{rep@theorem}}%
 {\end{rep@theorem}}}
\makeatother

\theoremstyle{plain}
\newtheorem{theorem}{Theorem}[section]
\newtheorem{lemma}[theorem]{Lemma}
\newtheorem{proposition}[theorem]{Proposition}

\newtheorem{conjecture}[theorem]{Conjecture}

\newtheorem{questions}[theorem]{Questions}

\theoremstyle{definition}
\newtheorem{definition}[theorem]{Definition}
\newtheorem{example}[theorem]{Example}

\theoremstyle{remark}
\newtheorem{remark}[theorem]{Remark}

\newreptheorem{theorem}{Theorem}

\newcommand{\CC}{\mathcal{C}}

\newcommand{\gin}{{\rm gin}}
\newcommand{\hgt}{{\rm ht}}

\newcommand{\depth}{{\rm depth}}

\newcommand{\cone}{{\rm cone}}
\newcommand{\link}{{\rm link}}
\newcommand{\mc}{\mathcal}

\newcommand{\mbf}{\mathbf}

\title{A combinatorial characterization of $S_2$ binomial edge ideals}

\author[D. Bolognini, A. Macchia, G. Rinaldo, F. Strazzanti]{Davide Bolognini, Antonio Macchia, Giancarlo Rinaldo, Francesco Strazzanti}

\address{{\small Davide Bolognini, Dipartimento di Ingegneria Industriale e Scienze Matematiche, Universit\`a Politecnica delle Marche, Via Brecce Bianche, 60131 Ancona, Italy}}
\email{{\small d.bolognini@univpm.it}}

\address{{\small Antonio Macchia, Fachbereich Mathematik und Informatik, Freie Universit\"at Berlin, Arnimallee 2, 14195 Berlin, Germany}}
\email{{\small macchia.antonello@gmail.com}}

\address{{\small Giancarlo Rinaldo, Dipartimento di Matematica e Informatica, Fisica e Scienze della Terra, Università di Messina, Viale Ferdinando Stagno d’Alcontres 31, 98166 Messina, Italy}}
\email{{\small giancarlo.rinaldo@unime.it}}

\address{{\small Francesco Strazzanti, Dipartimento di Matematica e Informatica, Fisica e Scienze della Terra, Università di Messina, Viale Ferdinando Stagno d’Alcontres 31, 98166 Messina, Italy}}
\email{{\small francesco.strazzanti@unime.it}}



\begin{document}

\begin{abstract}
Several algebraic properties of a binomial edge ideal $J_G$ can be interpreted in terms of combinatorial properties of its associated graph $G$. In particular, the so-called \textit{cut sets} of a graph $G$, special sets of vertices that disconnect $G$, play an important role since they are in bijection with the minimal prime ideals of $J_G$. In this paper we establish the first graph-theoretical characterization of binomial edge ideals $J_G$ satisfying Serre's condition $(S_2)$ by proving that this is equivalent to having $G$ \textit{accessible}, which means that $J_G$ is unmixed and the cut sets of $G$ form an accessible set system. The proof relies on the combinatorial structure of the Stanley--Reisner simplicial complex of a multigraded generic initial ideal of $J_G$, whose facets can be described in terms of cut sets. Another key step in the proof consists in proving the equivalence between accessibility and strong accessibility for the collection of cut sets of $G$ with $J_G$ unmixed. This result, interesting on its own, provides the first relevant class of set systems for which the previous two notions are equivalent.
\end{abstract}

\maketitle

\section{Introduction}

\textit{Accessible set systems} on a finite ground set $[n] = \{1,\dots,n\}$ are collections $\CC$ of subsets of $[n]$ such that every nonempty $S \in \CC$ contains an element $v$ for which $S \setminus \{v\} \in \CC$. They were introduced by Korte and Lovász in \cite{KL83} to study the structure of greedoids, a generalization of matroids. Moreover, accessible set systems further generalize greedoids producing a hierarchy of properties weaker than matroids as described in \cite{BHPW10}.

A relevant class of accessible set systems arises in the study of Cohen-Macaulay binomial edge ideals of graphs (see \cite{BMS18, BMS22}) and consists of the so-called \textit{cut sets} or \textit{cut-point sets} of a graph, which are sets of vertices that disconnect the graph in a ``minimal'' way. More precisely, if $G$ is a graph on the vertex set $[n]$, a set $S \subset [n]$ is a \textit{cut set} of $G$ if $c_G(S) > c_G(S \setminus \{i\})$ for every $i \in S$, where $c_G(S)$ denotes the number of connected components of $G - S$. We denote by $\CC(G)$ the collection of cut sets of $G$. Cut sets were first defined in \cite{HHHKR10} and \cite{O11} in order to combinatorially describe the minimal primes of binomial edge ideals, a class of binomial ideals $J_G \subset R = K[x_i, y_i : i \in [n]]$ of a polynomial ring $R$ generated by the binomials $x_iy_j - x_jy_i$ where $\{i, j\}$ are the edges of a finite simple graph $G$. A deep understanding of the combinatorial properties of cut sets turned out to be especially important in the study of unmixed and Cohen-Macaulay binomial edge ideals, see e.g. \cite{BMRS22, BMS18, BMS22, KS15, LMRR23, RR14, R19, SS22b}.

Notice that cut sets are special \textit{vertex-cuts} (see \cite[p. 230]{GYA19}) and have their own interest even outside of this algebraic context. Hence, they would deserve to be further studied in Graph Theory; however, we could not find any reference to the above definition of cut-point sets in the literature.

An important algebraic property of ideals that is often hard to grasp is Cohen-Macaulayness, hence it is interesting to find alternative ways of characterizing this property using a different language. In the case of combinatorially defined ideals, such as binomial edge ideals, it is natural to look for a characterization in terms of the underlying combinatorial object. In this direction, in \cite{BMS18} the first, second and fourth authors classified the bipartite graphs whose binomial edge ideals are Cohen-Macaulay, and a fundamental step in the proof is to show that the collection of cut sets of such graphs is an accessible set system. In a second paper, it turned out that the cut sets of Cohen-Macaulay binomial edge ideals of any graph form an accessible system \cite[Theorem 3.5]{BMS22}. Assuming $J_G$ unmixed, the converse holds for chordal and traceable graphs, see \cite[Theorems 6.4 and 6.8]{BMS22}. Supported by this evidence, they stated the following:

\begin{conjecture}[{\cite[Conjecture 1.1]{BMS22}}] \label{Conj.CM_accessible}
Let $G$ be a graph. Then, $J_G$ is Cohen–Macaulay if and only if $G$ is accessible.
\end{conjecture}

Here a graph $G$ is called \textit{accessible} if $\CC(G)$ is accessible and $J_G$ is unmixed, where the latter means that the minimal prime ideals of $J_G$ have the same height. For binomial edge ideals, the unmixedness is equivalent to $c_G(S) = |S| + c$ for every $S \in \CC(G)$, see \cite[Lemma 2.5]{RR14}, where $c$ is the number of connected components of $G$; hence, unmixedness is a combinatorial condition.

Conjecture \ref{Conj.CM_accessible} also holds for other classes of graphs described in \cite{LMRR23, SS22} and for graphs with a small number of vertices \cite[Theorem 1.2]{BMRS22}, providing further computational evidence. In \cite{LMRR23} Lerda, Mascia, Rinaldo and Romeo approached the conjecture by proving that the graph of a binomial edge ideal satisfying Serre's condition $(S_2)$ is accessible, strengthening \cite[Theorem 3.5]{BMS22}. At the same time they conjectured that the converse also holds, formulating a weaker version of Conjecture \ref{Conj.CM_accessible}:

\begin{conjecture}[{\cite[Conjecture 1]{LMRR23}}] \label{Conj.S2_accessible}
Let $G$ be a graph. Then, $R/J_G$ satisfies Serre's condition $(S_2)$ if and only if $G$ is accessible.
\end{conjecture}

Putting together the results of \cite{BMS22} and \cite{LMRR23}, we know that the following implications hold for all graphs:
\begin{equation}\label{Eq.implications}
\tag{$\ast$} R/J_G \text{ is Cohen-Macaulay } \Rightarrow R/J_G \text{ satisfies Serre's condition } (S_2) \Rightarrow G \text{ is accessible}.
\end{equation}

In this paper we prove the converse of the last implication, settling Conjecture \ref{Conj.S2_accessible} and providing the first characterization of $(S_2)$ binomial edge ideals in terms of the associated graph. This is a further step towards the original Conjecture \ref{Conj.CM_accessible}.

In Section \ref{S.MultiginComplex} we consider the main tool in our proof, the Stanley-Reisner simplicial complex $\Delta_G$ of the $\mathbb Z^n$-graded generic initial ideal $\gin(J_G)$ studied by Conca, De Negri and Gorla in \cite{CDG18}. In Lemma \ref{facets}, we explicitly describe the facets of the complex $\Delta_G$ and then show that any accessible graph $G$ has at least two free vertices, implying that $\Delta_G$ is an iterated cone, see Proposition \ref{P.gin_cone}.

In Section \ref{S.StronglyAccessible} we consider \textit{strongly accessible} set systems, first defined in \cite{BHPW10} as a subclass of accessible set systems. More precisely, these are set systems $\CC$ such that for every $S$, $T\in \CC$ with $S\subset T$ there exists $t\in T\setminus S$ such that $T\setminus\{t\}\in \CC$. The main result of this section is the following:

\begin{reptheorem}{T.stronglyAccessible}
Let $G$ be a graph with $J_G$ unmixed. Then $\CC(G)$ is accessible if and only if it is strongly accessible.
\end{reptheorem}

This result, which is interesting on its own, will play a key role in the proof of Conjecture \ref{Conj.S2_accessible} and it is the first relevant class in the literature for which these two properties are equivalent. Theorem \ref{T.stronglyAccessible} does not hold if $J_G$ is not unmixed, see Example \ref{E.accessibleNonStronglyAccessible}.

After some preliminary technical lemmas, Section \ref{S.S2} is dedicated to the proof of Conjecture \ref{Conj.S2_accessible}:

\begin{reptheorem}{T.S2BinomialEdgeIdeals}
Let $G$ be a graph. Then, $R/J_G$ satisfies Serre's condition $(S_2)$ if and only if $G$ is accessible.
\end{reptheorem}

As a consequence, we can rewrite the sequence of implications in \eqref{Eq.implications} as follows:
\begin{equation*}\label{Eq.implications2}
J_G \text{ Cohen-Macaulay } \Rightarrow R/J_G \text{ satisfies Serre's condition } (S_2) \Leftrightarrow G \text{ accessible}.
\end{equation*}

The study of the simplicial complex $\Delta_G$ is crucial in the proof of Theorem \ref{T.characterizationS2}. This suggests another approach to prove Conjecture \ref{Conj.CM_accessible}: it is enough to show that for every accessible graph $G$, the complex $\Delta_G$ is shellable, which in turns implies the Cohen-Macaulayness of $\gin(J_G)$ and hence of $J_G$.

\section{The simplicial complex of the \texorpdfstring{$\mathbb Z^n$}{Zn}-graded generic initial ideal of \texorpdfstring{$J_G$}{JG}} \label{S.MultiginComplex}

Let $G$ be a connected finite simple graph on the vertex set $[n]$ and let us use the notation defined above. Suppose that $S \subset [n]$ and let $G_1,\dots, G_{c_G(S)}$ be the connected components of $G - S$. Recall that $S$ is a \textit{cut set} of $G$ if $S = \emptyset$ or $c_G(S) > c_G(S \setminus \{i\})$ for every $i \in S$. Moreover, by \cite[Lemma 2.5]{RR14} $J_G$ is unmixed if and only if $c_G(S) = |S| + 1$ for every $S \in \CC(G)$.

We call a \textit{transversal} of $G - S$ a subset of vertices $W \subset V(G - S)$ such that $|W \cap G_i|=1$ for every $i = 1, \dots, c_G(S)$. We denote by $\mc T_G(S)$, or simply $\mc T(S)$, the set of all transversals of $G - S$.

The ideal $J_G$ is $\mathbb{Z}^n$-multigraded considering the natural $\mathbb{Z}^n$-graded structure on $R$ induced by\break $\deg(x_i)=\deg(y_i)=e_i \in \mathbb Z^n$ (where $(e_i)_j = 1$ if $j=i$ and $0$ otherwise) for $i = 1, \dots, n$. In \cite[Theorem 3.1]{CDG18}, Conca, De Negri and Gorla described the \textit{$\mathbb Z^n$-graded generic initial ideal} of $J_G$, denoted by $\gin(J_G)$, proving that $\gin(J_G)$ is a radical monomial ideal. Thus, $\gin(J_G)=I_{\Delta_G}$ is the Stanley-Reisner ideal of a simplicial complex $\Delta_G$ on the vertices $\{x_1,\dots,x_n,y_1,\dots,y_n\}$.

For each $S \in \mc C(G)$ and $W \in \mc T(S)$, we define the monomial
\[
F(S,W) = \prod_{i \notin S}y_i \prod_{j \in W}x_j \in R = K[x_1, \dots, x_n, y_1, \dots, y_n].
\]
Moreover, in what follows we confuse the set of vertices $\{y_i : i \notin S\} \cup \{x_j : j \in W\}$ of $\Delta_G$ with the monomial $F(S,W) \in R$.

\begin{lemma}\label{facets}
Let $G$ be a connected graph on the vertex set $[n]$. Then the facets of $\Delta_G$ are given by
\[
\{F(S,W) : S \in \mc C(G), W \in \mc T(S)\}.
\]
Moreover, $J_G$ is unmixed if and only if $\Delta_G$ is a pure simplicial complex of dimension $n$.
\end{lemma}

\begin{proof}
By \cite[Proposition 3.2]{CDG18}, the minimal primes of $\gin(J_G)$ are of the form $(x_i: i \notin W)+(y_j: j \notin A)$, where $A \subset [n]$ has the property that for every $i \in [n] \setminus A$, $c_G([n] \setminus (A \cup \{i\})) < c_G([n] \setminus A)$ and $W \subset A$ contains exactly one vertex in each connected component of the graph induced by $G$ on $A$. Setting $S = [n] \setminus A$, it is clear that the conditions above are equivalent to $S \in \mc C(G)$ and $W \in \mc T(S)$. Then the minimal primes of $\gin(J_G)$ are of the form $(x_i: i \notin W)+(y_j: y_j \in S)$, where $S \in \mc C(G), W \in \mc T(S)$. By \cite[Lemma 1.5.4]{HH10}, the facets of $\Delta_G$ are $\{F(S,W) : S \in \mc C(G), W \in \mc T(S)\}$, as desired.

We know that $J_G$ is unmixed if and only if $c_G(S)=|S|+1$ for every $S \in \mc C(G)$. This is equivalent to $|W|=|S|+1$, for every $W \in \mc T(S)$, i.e., $|F(S,W)| = n - |S| + |W| = n + 1$.
\end{proof}

\begin{example}
Let $G$ be the graph in Figure \ref{F.graphMultigin}, whose cut sets are $\CC(G)=\{\emptyset, \{2\}, \{1,2\}\}$. Notice that the sets of transversals corresponding to all cut sets of $G$ are: $\mc T(\emptyset) = \{\{1\}, \{2\}, \{3\}, \{4\}, \{5\}\}$, $\mc T(\{2\}) = \{\{1,3\}, \{3,4\}, \{3,5\}\}$ and $\mc T(\{1,2\}) = \{\{3,4,5\}\}$. By Lemma \ref{facets}, the facets of the Stanley-Reisner complex $\Delta_G$ of the ideal $\gin(J_G)$ are:
\begin{equation*}
\begin{aligned}
F(\emptyset, \{1\}) = x_1y_1y_2y_3y_4y_5,\\
F(\emptyset, \{2\}) = x_2y_1y_2y_3y_4y_5,\\
F(\emptyset, \{3\}) = x_3y_1y_2y_3y_4y_5,\\
F(\emptyset, \{4\}) = x_4y_1y_2y_3y_4y_5,\\
F(\emptyset, \{5\}) = x_5y_1y_2y_3y_4y_5,
\end{aligned}
\qquad
\begin{aligned}
F(\{2\}, \{1,3\}) = x_1x_3y_1y_3y_4y_5,\\
F(\{2\}, \{3,4\}) = x_3x_4y_1y_3y_4y_5,\\
F(\{2\}, \{3,5\}) = x_3x_5y_1y_3y_4y_5,\\
F(\{1,2\}, \{3,4,5\}) = x_3x_4x_5y_3y_4y_5.
\end{aligned}
\end{equation*}
Hence, $\Delta_G$ is pure of dimension $5$.

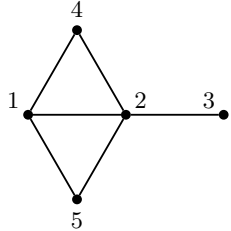
\begin{figure}[ht!]
\centering
\begin{tikzpicture}[scale=0.65]
\node[label={above left:{\small $1$}}] (a) at (0,0) {};
\node[label={above right:{\small $2$}}] (b) at (2,0) {};
\node[label={above left:{\small $3$}}] (c) at (4,0) {};
\node[label={above:{\small $4$}}] (d) at (1,1.7320508075688776) {};
\node[label={below:{\small $5$}}] (e) at (1,-1.7320508075688776) {};
\draw (c) -- (b) -- (a) -- (d) -- (b) -- (e) -- (a);
\end{tikzpicture}
\caption{A graph $G$}\label{F.graphMultigin}
\end{figure}
\end{example}

Following \cite{BMS22}, recall that a graph $G$ is \textit{accessible} if $J_G$ is unmixed and $\mc C(G)$ is an \textit{accessible set system}, i.e., for every non-empty $S \in \mc C(G)$ there exists $s \in S$ such that $S \setminus \{s\} \in \mc C(G)$. Moreover, a \textit{free vertex} of $G$ is a vertex that belongs to a unique maximal induced complete subgraph of $G$. As a first simple property of the complex $\Delta_G$ for $G$ accessible, we show that $\Delta_G$ is an iterated cone.

\begin{proposition}\label{P.gin_cone}
If $G$ is a connected accessible graph, then $G$ has at least two free vertices. In particular, $\Delta_G$ is a cone from $y_v$ for every free vertex $v$ of $G$.
\end{proposition}

\begin{proof}
Recall that a vertex $v$ of a graph $G$ is a \textit{cut vertex} of $G$ if $\{v\}$ is a cut set of $G$, and a graph is a \textit{block} if it does not have cut vertices. We can uniquely decompose $G$ into blocks, such that each two blocks of $G$ share at most one cut vertex of $G$. We say that a block is a \textit{terminal block} of $G$ if it contains exactly one cut vertex through which it is connected to other blocks of $G$.

If $G$ consists of a single block, then $G$ is complete by \cite[Remark 4.2]{BMS22} and all its vertices are free. Suppose that $G$ has at least two blocks, then it has at least a cut vertex; hence, it has at least two terminal blocks. We proceed by induction on the number $n \geq 3$ of vertices of $G$. If $n = 3$, then $G$ is a path of length two and it has two free vertices, one for each terminal block. Suppose $n > 3$ and consider a terminal block $B$ of $G$, whose only cut vertex is $v$. Then, if $x$ is a new vertex, the graph $H$ obtained by adding the edge $\{v, x\}$ to the graph $B$ is accessible by \cite[Proposition 3]{LMRR23}. Thus \cite[Lemma 4.9]{BMS22} yields that $H$ is a cone. More precisely, $H = \cone(v, \{x\} \sqcup (B - \{v\}))$, where $B - \{v\}$ is accessible by \cite[Theorem 4.8 (3)]{BMS22}. Since $B - \{v\}$ has fewer vertices than $G$, by induction it has a free vertex $w$. Clearly, $w$ is also a free vertex in $B$ (and then in $G$) because cones preserve free vertices. The last part of the statement follows by the definition of $\Delta_G$: if $v$ is a free vertex of $G$, then it is not contained in any cut set of $G$ by \cite[Lemma 2.1]{RR14}, hence $y_v$ belongs to all facets $F(S,W) = \prod_{i \notin S}y_i \prod_{j \in W}x_j$ of $\Delta_G$.
\end{proof}

When $G$ is bipartite, by \cite[Proposition 2.3]{BMS18} $G$ has exactly two free vertices requiring only $J_G$ unmixed.

\section{Accessible and strongly accessible set systems} \label{S.StronglyAccessible}

While looking for a combinatorial description of Cohen-Macaulay binomial edge ideals of bipartite graphs, in \cite{BMS18} we stumbled upon a characterization of the underlying graphs purely in terms of cut sets: condition (d) of \cite[Theorem 6.1]{BMS18} requires that $J_G$ is unmixed and that \textit{for every $S \in \CC(G)$ there exists $s \in S$ such that $S \setminus \{s\} \in \CC(G)$}, where the latter means that $\CC(G)$ is an \textit{accessible set system}. We later called these graphs \textit{accessible} and realized that these conditions seemed to be common to the graphs of all Cohen-Macaulay binomial edge ideals, leading to the formulation and the study of Conjecture \ref{Conj.CM_accessible}.

The notion of accessible set systems was originally defined by Korte and Lovász in \cite[p. 360]{KL83} as a superclass of \textit{greedoids} in an effort to further generalize matroids. More recently, Boley et al. \cite{BHPW10} refined even more this hierarchy, defining strongly accessible set systems as an intermediate notion between greedoids and accessible set systems.

\begin{definition}
A set system $\CC$ is called \textit{strongly accessible} if for every $S,T \in \CC$ with $S \subset T$, there exists $v \in T \setminus S$ such that $S \cup \{v\} \in \CC$.
\end{definition}

\begin{remark} \label{R.strongly_accessible}
The following properties are equivalent:
\begin{enumerate}
\item $\CC$ is strongly accessible;
\item for every $S,T \in \CC$ with $S \subset T$, there exists $t \in T \setminus S$ such that $T \setminus \{t\} \in \CC$;
\item for every $S,T \in \CC$ with $S \subset T$, there exists an ordering on $T \setminus S=\{t_1, \dots, t_r\}$ such that $S \cup \{t_1, \dots, t_i\} \in \CC$ for every $i=1, \dots, r$.
\end{enumerate}
Depending on the situation we will use one of the above alternative characterizations of strong accessibility. From Condition (2) it easily follows that every strongly accessible set system $\CC$ is accessible by choosing $S = \emptyset \in \CC$.
\end{remark}

\begin{example}\label{E.accessibleNonStronglyAccessible}
In general, an accessible set system is not necessarily strongly accessible. Let $\CC = \{\emptyset, \{2\}, \{3\}, \{1,2\},\break \{1,2,3\}\}$. One can easily check that $\CC$ is accessible, but it is not strongly accessible since $\{3\} \subset \{1,2,3\}$, and $\{1,3\}, \{2,3\} \notin \CC$.

This is also the case when we consider the collection of cut sets of a graph. In fact, the set system $\CC$ is also the collection of cut sets of the graph $G$ in Figure \ref{F.accessibleNonStronglyAccessible}. In this case, $J_G$ is not unmixed because $c_G(\{3\}) = 3 > 2 = |\{3\}| + 1$.

\begin{figure}[ht!]
\centering
\begin{tikzpicture}[scale=0.65]
\node[label={above left:{\small $1$}}] (a) at (0,0) {};
\node[label={above right:{\small $2$}}] (b) at (2,0) {};
\node[label={above left:{\small $3$}}] (c) at (4,0) {};
\node[label={above:{\small $4$}}] (d) at (1,1.7320508075688776) {};
\node[label={below:{\small $5$}}] (e) at (1,-1.7320508075688776) {};
\node[label={above:{\small $6$}}] (f) at (5.414213562373095,1.414213562373095) {};
\node[label={below:{\small $7$}}] (g) at (5.414213562373095,-1.414213562373095) {};
\draw (f) -- (c) -- (b) -- (a) -- (d) -- (b) -- (e) -- (a)
(c) -- (g);
\end{tikzpicture}
\caption{A graph $G$ with $\CC(G)$ accessible and not strongly accessible}\label{F.accessibleNonStronglyAccessible}
\end{figure}
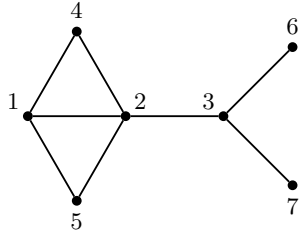
\end{example}

The aim of this section is to prove that accessibility and strong accessibility are indeed equivalent for the collection of cut sets of graphs whose binomial edge ideals are unmixed. This is a key result that will be used several times in the proof of Conjecture \ref{Conj.S2_accessible}.

In what follows, given a set $S \subset [n]$, we often say that a vertex $v \in S$ {\em reconnects} some connected components $G_1,\ldots, G_r$ of $G - S$ with $r \geq 2$, if adding back $v$ to $G - S$, together with all edges of $G$ incident in $v$, then $G_1,\ldots, G_r$ are in the same connected component. In particular, if $S$ is a cut set, any vertex in $S$ reconnects some connected components of $G - S$.

\begin{lemma} \label{L.CutsetContained}
Let $G$ be a graph on the vertex set $[n]$ and $U$ a proper subset of $[n]$. Then there exists a cut set $T$ of $G$ with $T \subset U$ such that $c_G(T) \geq c_G(U)$ and the following properties hold:
\begin{itemize}
\item[{\rm (i)}] if $W$ is part of a transversal of $G - U$, then it is also part of a transversal of $G - T$;
\item[{\rm (ii)}] if $u \in U \setminus T$, then $u$ does not reconnect any components of $G - U$;
\item[{\rm (iii)}] if $U$ is contained in a cut set, then $c_G(T) = c_G(U)$.
\end{itemize}
Moreover,
\begin{itemize}
\item[{\rm (iv)}] if $u_1, \dots, u_a \in U$ and $u_i$ does not reconnect any connected components in $G - (U \setminus \{u_1, \dots, u_{i-1}\})$ for every $i=1,\dots,a$, then we may assume $u_1,\dots,u_a \notin T$.
\end{itemize}
\end{lemma}

\begin{proof} Let $G_1, \dots, G_r$ be the connected components of $G - U$. We prove (i) and (iv) by induction on $|U| \geq 0$. Clearly, if $|U|=0$ the statement is trivial, so we assume $|U| > 0$. If $U$ is a cut set, it is enough to consider $T = U$ and the claim holds. Otherwise, there exists a vertex $u_1 \in U$ that does not reconnect any connected components in $G - U$, i.e., in $G - (U \setminus \{u_1\})$ either $u_1$ is an isolated vertex or it is adjacent only to vertices of some $G_j$ for exactly one $j$. Consider $U' = U \setminus \{u_1\}$. The connected components of $G - U'$ are either $G_1, \dots, G_r$ and the isolated vertex $\{u_1\}$ or the same components of $G - U$ except for $G_j$ that is replaced by the graph induced by $G$ on $V(G_j) \cup \{u_1\}$. It follows that $c_G(U') \geq c_G(U)$. By induction, there exists a cut set $T$ of $G$ contained in $U'$ (and then in $U$) such that $c_G(T) \geq c_G(U') \geq c_G(U)$. If $W$ is part of a transversal of $G - U$, then it is part of a transversal of $G - U'$ by construction and of $G - T$ by induction, hence (i) holds. Moreover, as for (iv), by induction we may assume that $u_2,\dots,u_a \notin T$; thus $u_1,\dots,u_a \notin T$.

Assume now that $u \in U \setminus T$ and by contradiction $u$ reconnects $G_i$ and $G_j$ in $G - U$ for some $i \neq j$. Let $w_i \in V(G_i)$ and $w_j \in V(G_j)$, then $\{w_i, w_j\}$ is part of a transversal of $G - U$, but it cannot be part of a transversal of $G - T$, against part (i). This proves (ii).

Finally, for (iii) we assume that $U$ is contained in $S \in \CC(G)$. If $u \in U \setminus T$, then it reconnects some connected components in $G-S$, thus in $G - (U \setminus \{u\})$ it is not an isolated vertex. This easily implies that $c_G(T) \leq c_G(U)$ because no new component is created by re-adding the vertices of $U \setminus T$ to $G - U$, hence the equality holds.
\end{proof}

\begin{example}
Note that the cut set $T$ found in Lemma \ref{L.CutsetContained} is not unique in general. For instance, let $G$ be the graph in Figure \ref{F.CutsetContained1} and $U = \{1,2,3,4,8\}$, which is not a cut set of $G$. There are two cut sets contained in $U$ and satisfying the claim of Lemma \ref{L.CutsetContained}: $\{2,4,8\}$ and $\{3,4,8\}$, see Figure \ref{F.CutsetContained}.

\begin{figure}[ht!]
\centering
\begin{tikzpicture}[scale=1.2]
\node[label={below:{\small $1$}}] (a) at (-0.5,0) {};
\node[label={below:{\small $2$}}] (b) at (1,0) {};
\node[label={below:{\small $3$}}] (c) at (2.5,0) {};
\node[label={above:{\small $4$}}] (d) at (0,0.75) {};
\node[label={above:{\small $5$}}] (e) at (1.25,0.75) {};
\node[label={above:{\small $6$}}] (f) at (-0.5,1.5) {};
\node[label={above:{\small $7$}}] (g) at (1,1.5) {};
\node[label={above:{\small $8$}}] (h) at (2.5,1.5) {};
\node[label={above:{\small $9$}}] (h) at (4,1.5) {};
\draw (-0.5,0) -- (1,0) -- (2.5,0) -- (1.25,0.75) -- (0,0.75) -- (1,0) -- (2.5,1.5) -- (1,1.5) -- (0,0.75)
(-0.5,1.5) -- (1,1.5)
(2.5,0) -- (2.5,1.5) -- (4,1.5)
(1.25,0.75) -- (2.5,1.5);
\end{tikzpicture}
\caption{A graph $G$ to illustrate Lemma \ref{L.CutsetContained}} \label{F.CutsetContained1}
\end{figure}
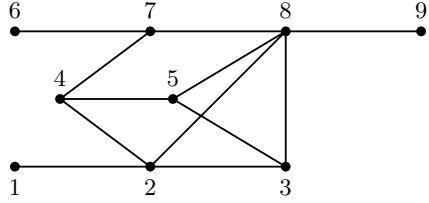

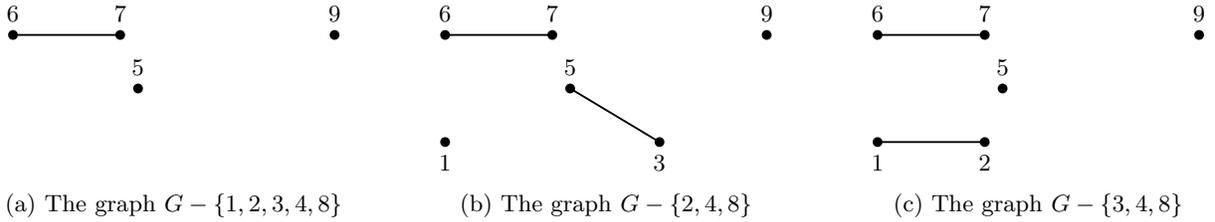
\begin{figure}[ht!]
\begin{subfigure}[c]{0.32\textwidth}
\centering
\begin{tikzpicture}[scale=0.95]
\phantom{\node[label={below:{\small $1$}}] (a) at (-0.5,0) {};}
\node[label={above:{\small $5$}}] (e) at (1.25,0.75) {};
\node[label={above:{\small $6$}}] (f) at (-0.5,1.5) {};
\node[label={above:{\small $7$}}] (g) at (1,1.5) {};
\node[label={above:{\small $9$}}] (h) at (4,1.5) {};
\draw (-0.5,1.5) -- (1,1.5);
\end{tikzpicture}
\caption{The graph $G - \{1,2,3,4,8\}$} \label{F.CutsetContained2}
\end{subfigure}
\begin{subfigure}[c]{0.32\textwidth}
\centering
\begin{tikzpicture}[scale=0.95]
\node[label={below:{\small $1$}}] (a) at (-0.5,0) {};
\node[label={below:{\small $3$}}] (c) at (2.5,0) {};
\node[label={above:{\small $5$}}] (e) at (1.25,0.75) {};
\node[label={above:{\small $6$}}] (f) at (-0.5,1.5) {};
\node[label={above:{\small $7$}}] (g) at (1,1.5) {};
\node[label={above:{\small $9$}}] (h) at (4,1.5) {};
\draw (2.5,0) -- (1.25,0.75)
(-0.5,1.5) -- (1,1.5);
\end{tikzpicture}
\caption{The graph $G - \{2,4,8\}$} \label{F.CutsetContained3}
\end{subfigure}
\begin{subfigure}[c]{0.32\textwidth}
\centering
\begin{tikzpicture}[scale=0.95]
\node[label={below:{\small $1$}}] (a) at (-0.5,0) {};
\node[label={below:{\small $2$}}] (b) at (1,0) {};
\node[label={above:{\small $5$}}] (e) at (1.25,0.75) {};
\node[label={above:{\small $6$}}] (f) at (-0.5,1.5) {};
\node[label={above:{\small $7$}}] (g) at (1,1.5) {};
\node[label={above:{\small $9$}}] (h) at (4,1.5) {};
\draw (-0.5,0) -- (1,0)
(-0.5,1.5) -- (1,1.5);
\end{tikzpicture}
\caption{The graph $G - \{3,4,8\}$} \label{F.CutsetContained4}
\end{subfigure}
\caption{An example of Lemma \ref{L.CutsetContained}} \label{F.CutsetContained}
\end{figure}
\end{example}

\begin{remark}
Lemma \ref{L.CutsetContained} (iv) implies that if $u \in S$ does not reconnect any connected components in $G - S$, then we may assume $u \notin T$. However, if $u_1, u_2 \in S$ do not reconnect any connected components in $G - S$, we cannot assume that both $u_1, u_2 \notin T$ without requiring that $u_2$ does not reconnect anything in $G - (S \setminus \{u_1\})$. For instance, consider the path $G$ with edge set $\{\{1,2\}, \{2,3\}, \{3,4\}\}$, and let $S = \{2,3\}$. Clearly, neither $2$ nor $3$ reconnect any components in $G - S$, but $T = S \setminus \{2,3\} = \emptyset$ does not satisfy the claim of Lemma \ref{L.CutsetContained} (this happens because $2$ reconnects two components of $G - (S \setminus \{3\})$ and vice versa exchanging the roles of $2$ and $3$). On the other hand, $T=\{2\}\subset S$   and $T=\{3\}\subset S$ satisfy the claim.
\end{remark}

Lemma \ref{L.CutsetContained} allows to give an immediate and completely different proof of the following lemma proved in \cite[Lemma 3.1]{BMRS22}, which we will frequently use in the rest of the paper.

\begin{lemma}\label{L.CutsetsInUnmixedGraphs}
Let $G$ be a connected graph with $J_G$ unmixed and let $S \subset V(G)$. Then the following properties are equivalent:
\begin{itemize}
\item[{\rm (1)}] $S$ is a cut set of $G$;
\item[{\rm (2)}] $c_{G}(S) = |S|+1$;
\item[{\rm (3)}] $c_{G}(S) \geq |S|+1$.
\end{itemize}
\end{lemma}

\begin{proof}
The implication $(1) \Rightarrow (2)$ follows by \cite[Lemma 2.5]{RR14} and $(2) \Rightarrow (3)$ is clear, hence we only need to prove that $(3)$ implies $(1)$. If $c_{G}(S) \geq |S|+1$, by Lemma \ref{L.CutsetContained} there exists a cut set $T \subset S$ such that $c_G(T) \geq c_G(S) \geq |S|+1 \geq |T|+1$. Therefore, the unmixedness of $J_G$ implies that the previous inequalities are all equalities and, in particular, $S=T \in \CC(G)$.
\end{proof}

\begin{remark} \label{R.union} Let $S$ be a cut set of a connected graph $G$ with $J_G$ unmixed.
\begin{enumerate}
\item If $s \in [n] \setminus S$, Lemma \ref{L.CutsetsInUnmixedGraphs} implies that $S \cup \{s\}$ is a cut set if and only if $s$ is a cut vertex of $G - S$.
\item If $s \in S$, Lemma \ref{L.CutsetsInUnmixedGraphs} implies that $S \setminus \{s\}$ is a cut set if and only if $s$ reconnects exactly two components of $G - S$. This was already shown in \cite[Lemma 4.14]{BMS22}.
\item Let $G_1, \dots, G_r$ be the connected components of $G - S$. There is at most one $s \in S$ that reconnects only $G_i$ and $G_j$ for two fixed $i$ and $j$. Indeed, if there are $s_1, \dots, s_a \in S$ reconnecting only $G_i$ and $G_j$ for $a \geq 2$, then $c_G(S \setminus \{s_1, \dots, s_a\})=c_G(S)-1=|S|>|S \setminus \{s_1, \dots, s_a\}|+1$, which contradicts the equivalence between $(2)$ and $(3)$ in Lemma \ref{L.CutsetsInUnmixedGraphs}.
\end{enumerate}
\end{remark}

\begin{theorem} \label{T.stronglyAccessible}
Let $G$ be a graph with $J_G$ unmixed. Then $\CC(G)$ is accessible if and only if it is strongly accessible.
\end{theorem}

\begin{proof}
First of all, notice that $\CC(G)$ is (strongly) accessible if and only if the same holds for each connected component of $G$. Therefore, it is enough to prove the claim for $G$ connected.

In Remark \ref{R.strongly_accessible}, we saw that every strongly accessible set system is accessible. Conversely, let $\CC(G)$ be accessible and suppose by contradiction that condition (2) in Remark \ref{R.strongly_accessible} does not hold. Let $T$ be a cut set of $G$ with the smallest cardinality for which there exists $S \in \CC(G)$ with $S\subset T$ and such that $T \setminus \{t\} \not\in \CC(G)$ for every $t \in T \setminus S$. Among all such cut sets $S$ we choose one with the largest cardinality.
Therefore, there are no cut sets between $S$ and $T$ and $|T \setminus S| \geq 2$.

Let $T\setminus S=\{t_1, \dots, t_r\}$, where $r\geq 2$, and let $G_1, \dots, G_{|T|+1}$ be the connected components of $G - T$. We may assume that the components of $G - S$ are the graphs induced by $G$ on each of the following sets of vertices:
\begin{gather*}
V(G_1) \cup \dots \cup V(G_{i_1})\cup \{t_1, \dots, t_{j_1}\}, \\
V(G_{i_1 +1}) \cup \dots \cup V(G_{i_2})\cup \{t_{j_1 +1}, \dots, t_{j_2}\}, \\
\vdots \\
V(G_{i_{a-1} +1}) \cup \dots \cup V(G_{i_a})\cup \{t_{j_{a-1} +1}, \dots, t_{j_a}\}, \\
V(G_{i_a+1}), \dots, V(G_{|T|+1}),
\end{gather*}
for some $a \geq 1$, $1 <i_1 < \dots < i_a \leq |T|+1$ and $j_1 < \dots < j_a=r$. Note that $i_a-a=r$ because\break $|T| + 1 - r = |S| + 1 = c_G(S) = a + (|T| + 1 - i_a)$.
We now prove that $a=1$. Setting $i_0=j_0=0$, we have
\[
\sum_{k=1}^a (i_k - i_{k-1}-1) = i_a - a=r=\sum_{k=1}^a (j_k - j_{k-1}),
\]
and then there exists an index $k$ such that $i_k - i_{k-1}-1 \geq j_k - j_{k-1}$. Let $U = S \cup \{t_{j_{k-1}+1}, \dots, t_{{j_{k}}}\}$. Then $c_G(U)=c_G(S)+(i_k-i_{k-1})-1 \geq c_G(S) + j_k - j_{k-1}=|S|+1+|\{t_{j_{k-1}+1}, \dots, t_{{j_{k}}}\}| = |U|+1$, and Lemma \ref{L.CutsetsInUnmixedGraphs} implies that $U$ is a cut set of $G$. Since $S \subsetneq U \subset T$, by the minimality of $T$ it follows that $U = T$, i.e., $a=1$ and $i_a = r+1$. This means that the connected components of $G - S$ are $G_{r+2}, \dots, G_{|T|+1}$ and the graph $H$ induced by $G$ on $V(G_1) \cup \dots \cup V(G_{r+1})\cup \{t_1, \dots, t_r\}$. As a consequence, every element of $S$ is adjacent to at least one vertex of a component among $G_{r+2},\dots, G_{|T|+1}$. Therefore, in $G - T$ every element of $S$ reconnects at least one component among $G_{r+2}, \dots, G_{|T|+1}$ to something else.

Clearly, $t_1$ reconnects at most $r+1$ components in $G - T$, say $G_1, \dots, G_{r'}$ with $r' \leq r + 1$. Consider the set
\[
Z=\{z \in T \mid z \text{ is not adjacent to any vertex in } V(G_{r'+1}) \cup \dots \cup V(G_{|T|+1})\}.
\]
Note that $t_1 \in Z$, and then by construction $T\setminus Z$ is a cut set strictly contained in $T$. The connected components of $G - (T \setminus Z)$ are $G_{r'+1}, \dots, G_{|T|+1}$ and the graph induced on $V(G_1) \cup \dots \cup V(G_{r'})\cup Z$, which is connected. Moreover, we proved above that every element of $S$ is not in $Z$, then $S \subset (T \setminus Z) \subsetneq T$. Hence, $T\setminus Z=S$ by the maximality of $S$ in $T$. It follows that $|T| - r + 1 = |S| + 1 = c_G(S) = c_G(T \setminus Z) = |T| + 1 - r' + 1$ which implies $r'=r+1$, i.e., $t_1$ reconnects $G_1, \dots, G_{r+1}$. Doing the same for $t_2, \dots, t_r$, we get that $t_i$ reconnects $G_1,\dots,G_{r+1}$ for every $i$.

Clearly $T$ is not empty, thus the accessibility of $G$ implies that there exists $t \in T$ such that $T \setminus \{t\} \in \CC(G)$; by assumption $t \in S$. This means that $t$ reconnects exactly two components in $G - T$ and at least two in $G - S$. We distinguish between two cases.

\textbf{Case 1.} First suppose that $t$ reconnects exactly two components also in $G - S$. In this case, $S \setminus \{t\} \in \CC(G)$ by Remark \ref{R.union} (2). By the minimality of $T$, there exists $v \in (T \setminus \{t\}) \setminus (S\setminus \{t\}) = T \setminus S$ such that $T \setminus \{t,v\} \in \CC(G)$. Since $v \in T \setminus S$, it reconnects $r+1$ components of $G - T$, and the components of $G - (T\setminus \{v\})$ are $G_{r+2}, \dots, G_{|T|+1}$ and the graph induced by $G$ on $V(G_1) \cup \dots \cup V(G_{r+1})\cup \{v\}$. Therefore, $t$ reconnects exactly two components of $G - (T\setminus \{v\})$ because $t \notin Z$. It follows that $c_G(T \setminus \{t,v\}) = |T|+1-r-1=|T|-r \leq |T| -2=|T \setminus \{t,v\}|$, and this contradicts the unmixedness of $J_G$.

\textbf{Case 2.} Suppose now that $t$ reconnects more than two components in $G - S$. This is possible only if it reconnects two components among $G_{r+2}, \dots, G_{|T|+1}$, say $G_{r+2}$ and $G_{r+3}$, and it is also adjacent to some vertices among $t_1, \dots, t_r$; then, $t$ reconnects three components in $G - S$: $G_{r+2}$, $G_{r+3}$, and $H$. In particular, in $G - T$ the vertex $t$ reconnects $G_{r+2}$ and $G_{r+3}$.

Notice that $S \setminus \{t\}$ is not a cut set because $c_G(S \setminus \{t\}) = c_G(S)-2 = |S|-1 = |S \setminus \{t\}|$, against the unmixedness of $J_G$. However, by Lemma \ref{L.CutsetContained} there exists a cut set $S'$ strictly contained in $S \setminus \{t\}$ and $c_G(S') \geq c_G(S \setminus \{t\}) = c_G(S) - 2 = |S| - 1 \geq (|S'|+2)-1 = |S'| + 1 = c_G(S')$, hence these are all equalities. This means that $S'=S \setminus \{t,u\}$ for some $u \in S$. Note that the connected components of $G - (S \setminus \{t\})$ are $G_{r+4}, \dots, G_{|T|+1}$ and the graph induced by $G$ on $V(G_1) \cup \dots \cup V(G_{r+3}) \cup \{t,t_1, \dots, t_r\}$. Moreover, by Lemma \ref{L.CutsetContained} (ii) $u$ does not reconnect any components of $G - (S \setminus \{t\})$; however, $u \in S$ reconnects some components among $G_{r+2}, G_{r+3}$, and $H$ in $G - S$ because $S$ is a cut set. On the other hand, since $t$ reconnects exactly $G_{r+2}$ and $G_{r+3}$ in $G - T$, Remark \ref{R.union} (3) implies that $u$ cannot reconnect only $G_{r+2}$ and $G_{r+3}$ in $G - T$; therefore, $u$ is adjacent to some vertices of $V(G_1) \cup \dots \cup V(G_{r+1})$ in $G$.

Now let us consider $S' = S \setminus \{t,u\} \subset T \setminus \{t\}$, where $S'$ and $T \setminus \{t\}$ are both cut sets of $G$ and $(T \setminus \{t\}) \setminus S' = \{t_1, \dots, t_r, u\}$. Note that, for any $i$, $T\setminus \{t,t_i\}$ cannot be a cut set because $c_G(T\setminus \{t,t_i\})=c_G(T)-1-r=|T|-r \leq |T|-2=|T\setminus \{t,t_i\}|$. Therefore, the minimality of $T$ implies that $T \setminus \{t,u\} \in \CC(G)$. Furthermore, considering $S' \subset T \setminus \{t,u\}$, again by the minimality of $T$, there exists $i$ such that $T'=T\setminus \{t,u,t_i\} \in \CC(G)$. On the other hand, the connected components of $G - T'$ are $G_{r+4}, \dots, G_{|T|+1}$ and the graph induced by $G$ on $V(G_1) \cup \dots \cup V(G_{r+3}) \cup \{t,u,t_i\}$, which is connected. Hence, $c_G(T')=c_G(T)-r-2=|T|-r-1 \leq |T|-3=|T'|$, which again contradicts the unmixedness of $J_G$.
\end{proof}

Even from a pure graph-theoretical perspective, it would be interesting to continue the combinatorial study of cut sets of graphs.

\begin{questions}\
\begin{itemize}
\item[{\rm (1)}] Which set systems $\CC \subset 2^{[n]}$ with $\emptyset \in \CC$ are the collections of cut sets of a graph on the vertex set $[n]$?
\item[{\rm (2)}] Which graphs $G$ are such that $\CC(G)$ is (strongly) accessible or a greedoid (even when $J_G$ is not unmixed)?
\item[{\rm (3)}] Without assuming $J_G$ unmixed, does the (strong) accessibility of $G$ translate any algebraic property of $J_G$?
\end{itemize}
\end{questions}

Concerning Question (1) we show that not all set systems containing $\emptyset$ are cut sets of a graph.

\begin{proposition}
There is no graph $G$ whose cut sets are $\CC(G) = \{\emptyset, \{1,2\}, \{1,3\}, \{2,3\}\}$.
\end{proposition}

\begin{proof}
First of all we notice that $1,2,3$ belong to the same connected component of $G$. In fact, if $1$ belongs to some component $H_1$ of $G$ and $2$ belongs to another component $H_2$ of $G$, then $1$ is a cut vertex of $H_1$ because $\{1,2\} \in \CC(G)$; therefore, $1$ is also a cut vertex of $G$, but this is a contradiction since $\{1\}\not\in \CC(G)$. Hence, it is enough to prove the claim for $G$ connected because the components not containing $1,2,3$ do not play any role.

Suppose that $c_G(\{1,2\}) = m \geq 2$ and let us call $G_1, \dots, G_m$ the components of $G - \{1,2\}$. If $1$ does not reconnect all $m$ components, then $2$ would be a cut vertex, a contradiction. Hence, both $1$ and $2$ reconnect $G_1, \dots, G_m$. Assume that $3 \in G_1$.

First suppose that $G_1 - \{3\}$ is connected. If both $1$ and $2$ are not adjacent to vertices of $G_1 - \{3\}$, then they are both adjacent to $3$. Thus, either $3$ is a cut vertex of $G$ or $G_1 = \{3\}$. The first case cannot occur by assumption. In the second case, $G - \{1,3\}$ is connected because $2$ is adjacent to $G_2, \dots, G_m$, a contradiction. Therefore we assume that $2$ is adjacent to some vertex in $G_1 - \{3\}$, but again $G - \{1,3\}$ is connected, a contradiction.

Hence, $G_1 - \{3\}$ is disconnected and let us call $H_1, \dots, H_r$ its connected components. Since $\{1,2,3\}$ is not a cut set and $3$ reconnects $H_1, \dots, H_r$ in $G - \{1,2,3\}$, we may assume that $1$ does not reconnect any component in $G - \{1,2,3\}$. However, since $1$ reconnects $G_2, \dots, G_m$, it follows that $m=2$ and $1$ does not reconnect any of $H_1, \dots, H_r$. Then $1$ does not reconnect any component in $G - \{1,3\}$, a contradiction.
\end{proof}

\section{Binomial edge ideals satisfying Serre's condition \texorpdfstring{$(S_2)$}{S2}} \label{S.S2}

Recall that a ring $A$ satisfies \textit{Serre's condition} $(S_r)$ if for every prime ideal $\mathfrak p$ of $A$, $\depth\, A_{\mathfrak p} \geq \min\{r, \hgt\, \mathfrak p\}$. Moreover, $A$ is \textit{Cohen-Macaulay} if it satisfies $(S_r)$ for every $r \geq 1$. Serre's conditions, especially $(S_2)$, are important in Commutative Algebra and Algebraic Geometry, see for instance \cite{H62} and \cite[p. 183]{M89}. Some results on $(S_2)$ binomial edge ideals can be found in \cite{BN17, LMRR23}. In particular, in \cite[Theorem 2]{LMRR23} it is proved that if $R/J_G$ satisfies Serre's condition $(S_2)$, then $G$ is accessible and in Conjecture \ref{Conj.S2_accessible} the reverse implication is conjectured.

In this section we settle the above conjecture by exploiting the structure of the simplicial complex $\Delta_G$. In particular, leveraging the following characterization of $(S_2)$ Stanley-Reisner rings, we prove that if $G$ is an accessible graph, then the links of the faces of $\Delta_G$ are either connected or zero-dimensional. Recall that, given a simplicial complex $\Delta$ on the vertex set $[n]$, its \textit{Stanley-Reisner ring} is $R/I_\Delta$, where $I_\Delta \subset R = K[x_1,\dots,x_n]$ is the ideal generated by the monomials corresponding to the minimal non-faces of $\Delta$. The dimension of $\Delta$ is\break $\dim(\Delta) = \max\{|F|-1 : F \in \Delta\}$ and the \textit{link} of a face $F$ of $\Delta$ is the simplicial complex
\[
\link_\Delta(F) = \{G \in \Delta : G \cap F = \emptyset, G \cup F \in \Delta \}.
\]
The next result follows from \cite[Theorem 1.4]{T07}, see also \cite[Corollary 2.4]{PSFTY14}.

\begin{theorem} \label{T.characterizationS2}
Let $\Delta$ be a simplicial complex. Then the following are equivalent:
\begin{itemize}
\item[{\rm (i)}] $R/I_\Delta$ satisfies Serre's condition $(S_2)$;
\item[{\rm (ii)}] $\link_\Delta(F)$ is connected for every face $F \in \Delta$ with $\dim(\link_\Delta(F)) \geq 1$.
\end{itemize}
\end{theorem}

Before proving Conjecture \ref{Conj.S2_accessible}, we need some technical lemmas that we illustrate through examples. These will help in proving certain cases in the proof of the main result. In the whole section we use the structure of the facets of $\Delta_G$ described in Lemma \ref{facets}.

\begin{example}
First we show how to connect two facets of $\link_{\Delta_G}(H)$ coming from two cut sets with transverals $W_1 \subset W_2$. Let $G$ be the accessible graph in Figure \ref{F.CutsetContained1} and set $T_1 = \{7,8\}, T_2 = \{3,4,8\} \in \CC(G)$. Notice that $T_1 \not \subset T_2$, $|T_1| < |T_2|$, and $W_1 = \{1,6,9\} \subset W_2 = \{1,5,6,9\}$ are transversals of $G - T_1$ and $G - T_2$ respectively, see Figure \ref{F.firstLemma}.

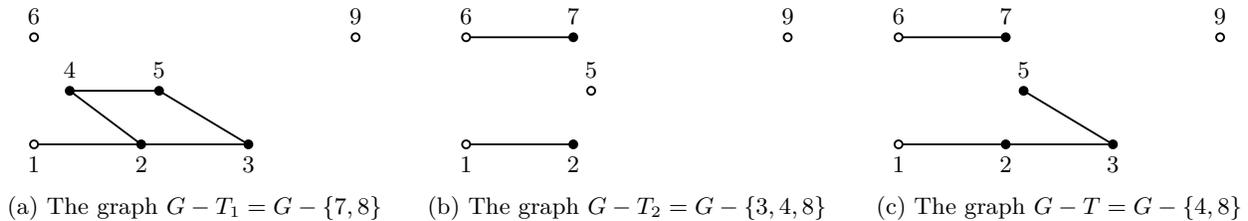
\begin{figure}[ht!]
\begin{subfigure}[c]{0.32\textwidth}
\centering
\begin{tikzpicture}[scale=0.95]
\draw (-0.5,0) -- (1,0) -- (2.5,0) -- (1.25,0.75) -- (0,0.75) -- (1,0);
\node[fill=white, label={below:{\small $1$}}] (a) at (-0.5,0) {};
\node[label={below:{\small $2$}}] (b) at (1,0) {};
\node[label={below:{\small $3$}}] (c) at (2.5,0) {};
\node[label={above:{\small $4$}}] (d) at (0,0.75) {};
\node[label={above:{\small $5$}}] (e) at (1.25,0.75) {};
\node[fill=white, label={above:{\small $6$}}] (f) at (-0.5,1.5) {};
\node[fill=white, label={above:{\small $9$}}] (h) at (4,1.5) {};
\end{tikzpicture}
\caption{The graph $G - T_1 = G - \{7,8\}$} \label{F.firstLemma1}
\end{subfigure}
\begin{subfigure}[c]{0.32\textwidth}
\centering
\begin{tikzpicture}[scale=0.95]
\draw (-0.5,0) -- (1,0)
(-0.5,1.5) -- (1,1.5);
\node[fill=white, label={below:{\small $1$}}] (a) at (-0.5,0) {};
\node[label={below:{\small $2$}}] (b) at (1,0) {};
\node[fill=white, label={above:{\small $5$}}] (e) at (1.25,0.75) {};
\node[fill=white, label={above:{\small $6$}}] (f) at (-0.5,1.5) {};
\node[label={above:{\small $7$}}] (g) at (1,1.5) {};
\node[fill=white, label={above:{\small $9$}}] (h) at (4,1.5) {};
\end{tikzpicture}
\caption{The graph $G - T_2 = G - \{3,4,8\}$} \label{F.firstLemma2}
\end{subfigure}
\begin{subfigure}[c]{0.32\textwidth}
\centering
\begin{tikzpicture}[scale=0.95]
\draw (-0.5,0) -- (1,0) -- (2.5,0) -- (1.25,0.75)
(-0.5,1.5) -- (1,1.5);
\node[fill=white, label={below:{\small $1$}}] (a) at (-0.5,0) {};
\node[label={below:{\small $2$}}] (b) at (1,0) {};
\node[label={below:{\small $3$}}] (c) at (2.5,0) {};
\node[label={above:{\small $5$}}] (e) at (1.25,0.75) {};
\node[fill=white, label={above:{\small $6$}}] (f) at (-0.5,1.5) {};
\node[label={above:{\small $7$}}] (g) at (1,1.5) {};
\node[fill=white, label={above:{\small $9$}}] (h) at (4,1.5) {};
\end{tikzpicture}
\caption{The graph $G - T = G - \{4,8\}$} \label{F.firstLemma3}
\end{subfigure}
\caption{An example to illustrate Lemma \ref{L.firstLemmaMainTheorem}} \label{F.firstLemma}
\end{figure}

In this case, the facets of $\Delta_G$ associated with these cut sets and transversals are
\begin{align*}
F_1 &= F(T_1, W_1) = (y_1 y_2 y_5 y_6 y_9 x_1 x_6 x_9) y_3 y_4 \text{ and }\\
F_2 &= F(T_2, W_2) = (y_1 y_2 y_5 y_6 y_9 x_1 x_6 x_9) y_7 x_5.
\end{align*}
Consider $H = F_1 \cap F_2 = y_1 y_2 y_5 y_6 y_9 x_1 x_6 x_9$ and we denote by $\widetilde F$ the face of $\link_{\Delta_G}(H)$ coming from a face $F$ of $\Delta_G$ with $H \subset F$. In particular, $\widetilde F_1 = y_3 y_4$ and $\widetilde F_2 = y_7 x_5$. Bearing in mind Theorem \ref{T.characterizationS2}, we want to find a facet $L = F(T,W)$ of $\Delta_G$ such that $\widetilde L \in \link_{\Delta_G}(H)$ and is connected both to $\widetilde F_1$ and to $\widetilde F_2$. Since $H \subset L$, we need $W_1 \subset W$ and $1,2,5,6,9 \notin T$, i.e., $T \subset T_1 \cup T_2$. Moreover, $\widetilde L$ is connected to $\widetilde F_1$ implies that $|W| = 3$ or $|W| = 4$, and one between $3$ and $4$ is not in $T$.

If $|W| = 4$, then $|T| = 3 = |T_1 \cup T_2| - 1$, and hence $T = (T_1 \cup T_2) \setminus \{3\} = T_1 \cup \{4\}$ or $T = (T_1 \cup T_2) \setminus \{4\} = T_1 \cup \{3\}$, but in both cases $T$ is not a cut set of $G$. Thus, $|W| = 3$ and we have that
\begin{itemize}
\item $\widetilde L$ is connected to $\widetilde F_1$ if $3 \notin T$ or $4 \notin T$, i.e., $T_2 \setminus T_1 \not \subset T$;
\item $\widetilde L$ is connected to $\widetilde F_2$ if $7 \notin T$, i.e., $T_1 \setminus T_2 \not \subset T$.
\end{itemize}

Indeed $T = \{4,8\}$ is a cut set of $G$ and $W = W_1$ is a transversal of $G - T$. The facet corresponding to $T$ and $W$ is $L = F(T, W) = (y_1 y_2 y_5 y_6 y_9 x_1 x_6 x_9) y_3 y_7$ and contains $H$, i.e., $\widetilde L$ is a facet of $\link_{\Delta_G}(H)$. We can see that in $\link_{\Delta_G}(H)$, $\widetilde L = y_3 y_7$ is connected to $\widetilde F_1$ through $y_3$ and to $\widetilde F_2$ through $y_7$:
\[
\widetilde F_1 = y_3 y_4\ \frac{\ \ y_3\ \ }{ }\ \widetilde L = y_3 y_7\ \frac{\ \ y_7\ \ }{ }\ \widetilde F_2 = y_7 x_5.
\]
\end{example}

\begin{lemma}\label{L.firstLemmaMainTheorem}
Let $G$ be a connected accessible graph and $T_1,T_2 \in \CC(G)$, with $T_1 \not \subset T_2$, $|T_1| < |T_2|$ such that $W_1 \subset W_2$ are transversals of $G - T_1$ and of $G - T_2$, respectively. Suppose that
\begin{itemize}
\item[{\rm (i)}] $T_1 \cup \{v\} \notin \CC(G)$ for every $v \in T_2 \setminus T_1$;
\item[{\rm (ii)}] $W_2 \setminus W_1 \not \subset T_1$.
\end{itemize}
Then there exists a cut set $T \subset T_1 \cup T_2$ of $G$ such that $T_1 \not \subset T$, $T_2 \setminus T_1 \not \subset T$ and with a transversal containing $W_1$. Moreover, given $f \in W_2 \setminus W_1$, $f \notin T_1$, there exists $c \in T_2 \setminus T_1$, $c \notin T$, such that $T \cup \{c\} \in \CC(G)$ and $W_1 \cup \{f\}$ is part of a transversal of $G - (T \cup \{c\})$.
\end{lemma}

\begin{proof}
First of all, notice that $\CC(G)$ is strongly accessible by Theorem \ref{T.stronglyAccessible}. Let $B_1,\dots,B_r$ be the connected components of $G - T_1$ and let $f \in W_2 \setminus W_1$ with $f \notin T_1$. Without loss of generality, we may assume that $f \in V(B_1)$ and let $b \in W_1 \cap V(B_1)$. Since $b,f \in W_2$, then $b,f \notin T_2$. It follows that $b,f \in V(B_1 - T_2)$. Since the connected components of $G - (T_1 \cup (T_2 \cap V(B_1)))$ are $B_2, \dots, B_r$ and the components of $B_1 - T_2$, where $B_1 - T_2$ is not connected. In fact, since $b,f \in W_2$, they are in different components of $G - T_2$, then cannot be in the same component of $B_1 - T_2$.

Let $C \subset T_2 \setminus T_1$ be a minimal subset (by inclusion) of $T_2 \cap V(B_1)$ for which $B_1 - C$ is not connected and $f,b$ are in different connected components. Notice that $|C| \geq 2$, since if $|C|=1$, then $c_G(T_1 \cup C) \geq c_G(T_1)+1 = |T_1|+2 = |T_1 \cup C|+1$ and so $T_1 \cup C$ would be a cut set by Lemma \ref{L.CutsetsInUnmixedGraphs}, against assumption (i). Moreover, $T_1 \cup C \notin \CC(G)$, because otherwise the strong accessibility of $\CC(G)$ would imply again that there exists $c \in C$ such that $T_1 \cup \{c\} \in \CC(G)$. By the minimality of $C$, it follows that every element of $C$ reconnects some components of $G - (T_1 \cup C)$. Since $T_1 \cup C \notin \CC(G)$, it means that there are some vertices $v_1, \dots, v_t \in T_1$, with $t \geq 1$, that do not reconnect anything in $G - (T_1 \cup C)$. However, $v_1, \dots, v_t \in T_1$ reconnect some components of $G - T_1$ to $B_1$ and in $B_1$ they are adjacent only to vertices in $C$, thus in $G - T_1$ every $v_i$ reconnects $B_1$ to exactly another connected component of $G - T_1$. Moreover, if $i \neq j$, then $v_i$ and $v_j$ do not reconnect the same component of $G - T_1$ to $B_1$, because there is at most one vertex of $T_1$ reconnecting two fixed components of $G - T_1$, see Remark \ref{R.union} (3). Thus $c_G(T_1 \setminus \{v_1, \dots, v_t\}) = c_G(T_1) - t = |T_1|+1-t = |T_1 \setminus \{v_1, \dots, v_t\}| + 1$, and hence $T_1 \setminus \{v_1, \dots, v_t\} \in \CC(G)$ by Lemma \ref{L.CutsetsInUnmixedGraphs}. Let $U = (T_1 \setminus \{v_1, \dots, v_t\})\cup C$. By construction, every vertex of $U$ reconnects at least two connected components of $G - U$, i.e., $U \in \CC(G)$ as well.

Since $T_1 \setminus \{v_1, \dots, v_t\} \subset U = (T_1 \setminus \{v_1, \dots, v_t\})\cup C$ are cut sets of $G$ and $\CC(G)$ is strongly accessible, by Remark \ref{R.strongly_accessible} (2) there exists $c \in C$ such that $T = U \setminus \{c\} = (T_1 \setminus \{v_1, \dots, v_t\})\cup (C \setminus \{c\})$ is a cut set of $G$. Clearly $T \subset U \subset T_1 \cup T_2$ and $T_1 \not \subset T$ because $t \geq 1$. Moreover, we know that $c \in C \subset T_2 \setminus T_1$ and $c \notin T$. In addition, since $v_1,\dots,v_t$ do not reconnect any components of $G - (T_1 \cup C)$, $W_1$ is part of a transversal of $G - T = G - (U \setminus \{c\})$. Finally, note that $T \cup \{c\} = U$ is a cut set and $W_1 \cup \{f\}$ is part of a transversal of $G - U$, since $f$ and $b$ are in different components of $G - U$ by definition of $C$.
\end{proof}

\begin{example}
Now we show how to connect two facets of $\link_{\Delta_G}(H)$ when the assumptions of Lemma \ref{L.firstLemmaMainTheorem} are not satisfied. Let $G$ be the accessible graph in Figure \ref{F.CutsetContained1} and set $T_1 = \{2,4,8\}, T_2 = \{2,5,7,8\} \in \CC(G)$. Notice that $T_1 \not \subset T_2$, and $W_1 = \{1,3,6,9\}$ is a transversal of $G - T_1$ contained in $W_2 = \{1,3,4,6,9\}$, which is a transversal of $G - T_2$, see Figure \ref{F.secondLemma1}.

\begin{figure}[ht!]
\begin{subfigure}[c]{0.32\textwidth}
\centering
\begin{tikzpicture}[scale=0.95]
\draw (2.5,0) -- (1.25,0.75)
(-0.5,1.5) -- (1,1.5);
\node[fill=white, label={below:{\small $1$}}] (a) at (-0.5,0) {};
\node[fill=white, label={below:{\small $3$}}] (c) at (2.5,0) {};
\node[label={above:{\small $5$}}] (e) at (1.25,0.75) {};
\node[fill=white, label={above:{\small $6$}}] (f) at (-0.5,1.5) {};
\node[label={above:{\small $7$}}] (g) at (1,1.5) {};
\node[fill=white, label={above:{\small $9$}}] (h) at (4,1.5) {};
\end{tikzpicture}
\caption{The graph $G - T_1 = G - \{2,4,8\}$} \label{F.secondLemma1}
\end{subfigure}
\begin{subfigure}[c]{0.32\textwidth}
\centering
\begin{tikzpicture}[scale=0.95]
\node[fill=white, label={below:{\small $1$}}] (a) at (-0.5,0) {};
\node[fill=white, label={below:{\small $3$}}] (c) at (2.5,0) {};
\node[fill=white, label={above:{\small $4$}}] (d) at (0,0.75) {};
\node[fill=white, label={above:{\small $6$}}] (f) at (-0.5,1.5) {};
\node[fill=white, label={above:{\small $9$}}] (h) at (4,1.5) {};
\end{tikzpicture}
\caption{The graph $G - T_2 = G - \{2,5,7,8\}$} \label{F.secondLemma2}
\end{subfigure}
\begin{subfigure}[c]{0.32\textwidth}
\centering
\begin{tikzpicture}[scale=0.95]
\draw (2.5,0) -- (1.25,0.75) -- (0,0.75);
\node[fill=white, label={below:{\small $1$}}] (a) at (-0.5,0) {};
\node[fill=white, label={below:{\small $3$}}] (c) at (2.5,0) {};
\node[label={above:{\small $4$}}] (d) at (0,0.75) {};
\node[label={above:{\small $5$}}] (e) at (1.25,0.75) {};
\node[fill=white, label={above:{\small $6$}}] (f) at (-0.5,1.5) {};
\node[fill=white, label={above:{\small $9$}}] (h) at (4,1.5) {};
\end{tikzpicture}
\caption{The graph $G - T = G - \{2,7,8\}$} \label{F.secondLemma3}
\end{subfigure}
\caption{An example to illustrate Lemma \ref{L.secondLemmaMainTheorem}} \label{F.secondLemma}
\end{figure}
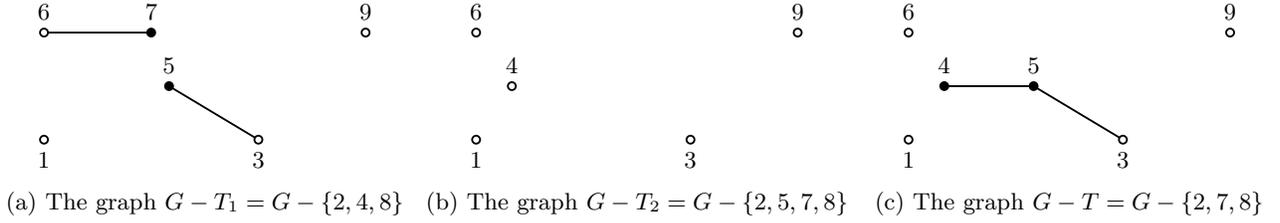

In this case, the facets of $\Delta_G$ associated with the cut sets and transversals above are
\begin{align*}
F_1 &= F(T_1, W_1) = (y_1 y_3 y_6 y_9 x_1 x_3 x_6 x_9) y_5 y_7 \text{ and }\\
F_2 &= F(T_2, W_2) = (y_1 y_3 y_6 y_9 x_1 x_3 x_6 x_9) y_4 x_4.
\end{align*}
Consider $H = F_1 \cap F_2 = y_1 y_3 y_6 y_9 x_1 x_3 x_6 x_9$, hence $\widetilde F_1 = y_5 y_7, \widetilde F_2 = y_4 x_4$ are facets of $\link_{\Delta_G}(H)$. In this case, there exists a cut set $T = \{2,7,8\}$ of $G$ and $W_1$ is a transversal of $G - T$. The facet corresponding to $T$ and $W_1$ is $L = F(T, W_1) = (y_1 y_3 y_6 y_9 x_1 x_3 x_6 x_9) y_4 y_5$ and contains $H$, i.e., $\widetilde L$ is a facet of $\link_{\Delta_G}(H)$. We can see that in $\link_{\Delta_G}(H)$, $\widetilde L = y_4 y_5$ is connected to $\widetilde F_1$ through $y_5$ and to $\widetilde F_2$ through $y_4$:
\[
\widetilde F_1 = y_5 y_7\ \frac{\ \ y_5\ \ }{ }\ \widetilde L = y_4 y_5\ \frac{\ \ y_4\ \ }{ }\ \widetilde F_2 = y_4 x_4.
\]
Note that, $f = 4 \in T_1 \setminus T_2$ is such that $N_G(f) = \{2,5,7\} \subset T_1 \cup T_2$. In the next lemma we will see that the existence of such $f$ guarantees the existence of the cut set $T$.
\end{example}

\begin{lemma}\label{L.secondLemmaMainTheorem}
Let $G$ be a connected accessible graph and $T_1,T_2 \in \CC(G)$, with $T_1 \not \subset T_2$ and consider two transversal $W_1$, $W_2$ of $G - T_1$ and $G - T_2$, respectively. Assume $W_1 \subset W_2$ and suppose that there exists $f \in T_1 \setminus T_2$ such that $N_G(f) \subset T_1 \cup T_2$. Then there exists a cut set $T \subset T_1 \cup T_2$ of $G$ such that $T_1 \setminus T_2 \not \subset T$, $T_2 \setminus T_1 \not \subset T$ and with a transversal containing $W_1$.
\end{lemma}

\begin{proof}
Since $T_1$ is a cut set and $f \in T_1$, $f$ reconnects some connected components of $G - T_1$, say $B_1, \dots, B_t$ for some $t \geq 2$. By our assumption, $N_{B_i}(f) \subset T_2$ for every $1 \leq i \leq t$. Consider the set $U = (T_1 \setminus \{f\}) \cup N_{B_1}(f) \cup \dots \cup N_{B_{t-1}}(f)$. By construction, $W_1$ is part of a transversal of $G - U$,  because $W_1 \cap N_{B_j}(f)= \emptyset$, for every $1 \leq j \leq t-1$, otherwise $\emptyset \neq W_1 \cap T_2 \subset W_2 \cap T_2$, a contradiction. By Lemma \ref{L.CutsetContained}, there exists a cut set $T \subset U$ of $G$ such that $W_1$ is part of a transversal of $G - T$. Notice that $f \in T_1 \setminus T_2$ and $f \notin T$, hence $T_1 \setminus T_2 \not \subset T$. Moreover, there exists $c \in N_{B_t}(f) \subset T_2$ and $c \notin T_1$, because $c \in B_t$, a connected component of $G - T_1$; on the other hand, $c \notin U \supset T$, then $T_2 \setminus T_1 \not \subset T$.
\end{proof}

\begin{example}
Finally we show how to connect two facets of $\link_{\Delta_G}(H)$ coming from two cut sets with a common transversal. Let $G$ be the accessible graph in Figure \ref{F.CutsetContained1} and set $T_1 = \{3,4,8\}, T_2 = \{2,7,8\} \in \CC(G)$. Notice that $|T_1 \setminus T_2| = |\{3,4\}| \geq 2$ and $W = \{1,5,6,9\}$ is a transversal of both $G - T_1$ and $G - T_2$, see Figure \ref{F.thirdLemma}.

\begin{figure}[ht!]
\begin{subfigure}[c]{0.32\textwidth}
\centering
\begin{tikzpicture}[scale=0.95]
\draw (-0.5,0) -- (1,0)
(-0.5,1.5) -- (1,1.5);
\node[fill=white, label={below:{\small $1$}}] (a) at (-0.5,0) {};
\node[label={below:{\small $2$}}] (b) at (1,0) {};
\node[fill=white, label={above:{\small $5$}}] (e) at (1.25,0.75) {};
\node[fill=white, label={above:{\small $6$}}] (f) at (-0.5,1.5) {};
\node[label={above:{\small $7$}}] (g) at (1,1.5) {};
\node[fill=white, label={above:{\small $9$}}] (h) at (4,1.5) {};
\end{tikzpicture}
\caption{The graph $G - T_1 = G - \{3,4,8\}$} \label{F.thirdLemma1}
\end{subfigure}
\begin{subfigure}[c]{0.32\textwidth}
\centering
\begin{tikzpicture}[scale=0.95]
\draw (2.5,0) -- (1.25,0.75) -- (0,0.75);
\node[fill=white, label={below:{\small $1$}}] (a) at (-0.5,0) {};
\node[label={below:{\small $3$}}] (c) at (2.5,0) {};
\node[label={above:{\small $4$}}] (d) at (0,0.75) {};
\node[fill=white, label={above:{\small $5$}}] (e) at (1.25,0.75) {};
\node[fill=white, label={above:{\small $6$}}] (f) at (-0.5,1.5) {};
\node[fill=white, label={above:{\small $9$}}] (h) at (4,1.5) {};
\end{tikzpicture}
\caption{The graph $G - T_2 = G - \{2,7,8\}$} \label{F.thirdLemma2}
\end{subfigure}
\begin{subfigure}[c]{0.32\textwidth}
\centering
\begin{tikzpicture}[scale=0.95]
\draw (2.5,0) -- (1.25,0.75)
(-0.5,1.5) -- (1,1.5);
\node[fill=white, label={below:{\small $1$}}] (a) at (-0.5,0) {};
\node[label={below:{\small $3$}}] (c) at (2.5,0) {};
\node[fill=white, label={above:{\small $5$}}] (e) at (1.25,0.75) {};
\node[fill=white, label={above:{\small $6$}}] (f) at (-0.5,1.5) {};
\node[label={above:{\small $7$}}] (g) at (1,1.5) {};
\node[fill=white, label={above:{\small $9$}}] (h) at (4,1.5) {};
\end{tikzpicture}
\caption{The graph $G - T = G - \{2,4,8\}$} \label{F.thirdLemma3}
\end{subfigure}
\caption{An example to illustrate Lemma \ref{L.thirdLemmaMainTheorem}} \label{F.thirdLemma}
\end{figure}
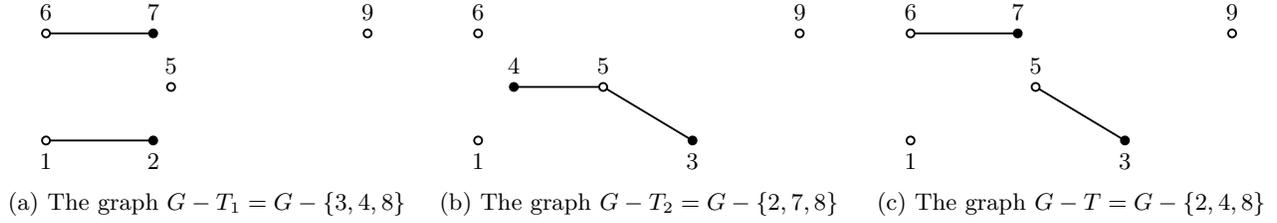

In this case, the facets of $\Delta_G$ associated to the cut sets and transversals above are
\begin{align*}
F_1 &= F(T_1, W_1) = (y_1 y_5 y_6 y_9 x_1 x_5 x_6 x_9) y_2 y_7 \text{ and }\\
F_2 &= F(T_2, W_2) = (y_1 y_5 y_6 y_9 x_1 x_5 x_6 x_9) y_3 y_4.
\end{align*}
Consider $H = F_1 \cap F_2 = y_1 y_5 y_6 y_9 x_1 x_5 x_6 x_9$, hence $\widetilde F_1 = y_2 y_7$ and $\widetilde F_2 = y_3 y_4$ are facets of $\link_{\Delta_G}(H)$. Notice that there is no cut set of $G$ contained in $T_1 \setminus T_2$ with a transversal strictly containing $W$. Therefore, to find a facet\break $L = F(T,X)$ such that $\widetilde L$ is connected both to $\widetilde F_1$ and $\widetilde F_2$ we need a cut set $T$ of $G$ such that\break $T \subset T_1 \cup T_2 = \{2,3,4,7,8\}$, satisfying $T_1 \setminus T_2 \not \subset T$, $T_2 \setminus T_1 \not \subset T$ and such that $X = W$. A cut set with these properties is $T = \{2,4,8\}$.

The facet $L = F(T, W) = (y_1 y_5 y_6 y_9 x_1 x_5 x_6 x_9) y_3 y_7$ contains $H$. We can see that in $\link_{\Delta_G}(H)$, $\widetilde L = y_3 y_7$ is connected to $\widetilde F_1$ through $y_7$ and to $\widetilde F_2$ through $y_3$:
\[
\widetilde F_1 = y_2 y_7\ \frac{\ \ y_7\ \ }{ }\ \widetilde L = y_3 y_7\ \frac{\ \ y_3\ \ }{ }\ \widetilde F_2 = y_3 y_4.
\]
\end{example}

\begin{lemma}\label{L.thirdLemmaMainTheorem}
Let $G$ be a connected accessible graph, $T_1,T_2 \in \CC(G)$, with $r = |T_1 \setminus T_2| \geq 2$, and $W$ a transversal of both $G - T_1$ and $G - T_2$. Suppose that there is no cut set of $G$ contained in $T_1 \cup T_2$ with a transversal strictly containing $W$. Then there exists a cut set $T \subset T_1 \cup T_2$ of $G$ satisfying $T_1 \setminus T_2 \not \subset T$, $T_2 \setminus T_1 \not \subset T$ and such that $W$ is a transversal of $G - T$.
\end{lemma}

\begin{proof}
Let $B_1, \dots, B_{|W|}$ be the connected components of $G - T_1$, where $|W| = c_G(T_1)$. We claim that $B_i - T_2$ is connected for $1 \leq i \leq |W|$. In fact, if there is an index $i$ for which $B_i - T_2$ is not connected, then $G - (T_1 \cup T_2)$ has more than $c_G(T_1)$ connected components and by Lemma \ref{L.CutsetContained} there exists a cut set $H \subseteq T_1 \cup T_2$ with a transversal of $G-H$ strictly containing $W$, against our assumption. The same argument works replacing $T_2$ with $T_1$.

We suppose first that there exists $e \in T_1 \setminus T_2$ such that $T_2 \setminus T_1 \not \subset N_G(e)$. Assume that $e$ reconnects $B_1,\ldots,B_i$. Consider $U=(T_1 \setminus \{e\}) \cup (N_{B_1}(e) \cap T_2) \cup \ldots \cup (N_{B_i}(e) \cap T_2)$ and let $T \subset U$ be a cut set obtained by Lemma \ref{L.CutsetContained}. Since $W$ is a transversal of $G - T_1$ and $G - T_2$, it is part of a transversal of $G - U$ and hence of $G - T$ by Lemma \ref{L.CutsetContained} (i). By our assumption, there is no cut set contained in $T_1 \cup T_2$ having a transversal strictly containing $W$, and thus, $W$ is a transversal of $G - T$. It is clear that $T \subset T_1 \cup T_2$ and $T_1 \setminus T_2 \not \subset T$. Since $T_2 \setminus T_1 \not \subset N_G(e)$, we have $T_2 \setminus T_1 \not \subset N_{B_1}(e) \cup \ldots \cup N_{B_i}(e)$; it follows that $T_2 \setminus T_1 \not \subset T$ and then $T$ is the cut set we were looking for.

Therefore, we may assume that there are no $e \in T_1 \setminus T_2$ such that $T_2 \setminus T_1 \not \subset N_G(e)$. This means that
\begin{equation}\label{Eq.adjacentVertices}
\text{if $e \in T_1 \setminus T_2$ and $c \in T_2 \setminus T_1$, then $e$ and $c$ are adjacent in $G$.}
\end{equation}

Let $T \subset T_1 \cap T_2$ be a cut set obtained by Lemma \ref{L.CutsetContained}. Let $T_2=T \cup S \cup (T_2 \setminus T_1)$, for some set $S$ with $|S| = s \geq 0$. By Theorem \ref{T.stronglyAccessible}, $\CC(G)$ is strongly accessible; since $T \subset T_2$ are cut sets, by Remark \ref{R.strongly_accessible} (3) there exists an ordering of the elements of $T_2 \setminus T = S \cup (T_2 \setminus T_1)=\{g_1,\ldots,g_{r+s}\}$ such that $T \cup \{g_1, \dots, g_i\} \in \CC(G)$, for every $0 \leq i \leq r+s$, where $r = |T_1 \setminus T_2| = |T_2 \setminus T_1|$ because $|T_1| = |T_2|$.

Suppose that $B_1,\ldots,B_a$ are all the connected components of $G - T_1$ such that $V(B_j) \cap T_2 \neq \emptyset$ for every $1 \leq j \leq a$. We claim that $|V(B_j) \cap T_2|=1$ for $1 \leq j \leq a$, i.e. $a=r$. In fact, if $c,d \in T_2 \setminus T_1$ are in the same component of $G \setminus T_1$, say $B_1$, then without loss of generality we have $c = g_h$ and $d = g_k$ with $h < k$. Thus $H=T \cup \{g_1, \dots, g_h\}$ is a cut set of $G$ not containing $d$. Therefore, the vertex $c$ does not reconnect any components of $G\setminus H$ because it is only adjacent to all the vertices of $T_1 \setminus T_2$ (by our assumption) and to some vertices of $B_1$, but those vertices belong to the same connected component of $G - H$, since $B_1 - T_2$ is connected, $d \in V(B_1)$ and it is also adjacent to all the vertices of $T_1 \setminus T_2$ since $T_1 \setminus T_2 \subset N_G(d)$ by assumption. This yields a contradiction, thus $a=r$. By symmetry,
\begin{equation}\label{Eq.e_jInComponentsOfG-T_2}
\text{the elements of $T_1 \setminus T_2$ belong to $r$ pairwise distinct connected components of $G - T_2$.}
\end{equation}

Let $T_2 \setminus T_1=\{c_1,\ldots,c_r\}$ and $T_1 \setminus T_2=\{e_1,\ldots,e_r\}$. Recall that by \eqref{Eq.adjacentVertices}, $c_i$ and $e_j$ are adjacent for every $i,j \in [r]$. Moreover, without loss of generality we assume that
\begin{equation*}\label{Eq.c_jInB_j}
c_j \in B_j \text{ for every } j = 1,\dots,r.
\end{equation*}

Therefore, in $G-T_1$ the vertex $e_1$ reconnects at least $B_1, \dots, B_r$. We now show that it is enough to assume $e_1$ also reconnects exactly one component of $G - T_1$ other than $B_1,\dots,B_r$. In fact, if in $G - T_1$ the vertex $e_1$ is adjacent to $v_i, v_j \notin T_2$, with $v_i \in V(B_i)$ and $v_j \in V(B_j)$ for $i \neq j$, then in $G-T_2$ there is a connected component containing $V(B_i - T_2) \cup V(B_j - T_2) \cup \{e_1\}$ and $W$ cannot be a transversal of $G - T_2$ because it would intersect this component in two vertices. Hence, $N_{G - T_1}(e_1) \setminus T_2$ is either empty or contained in a unique connected component $B_t$ of $G - T_1$. Assume first that $N_{G - T_1}(e_1) \setminus T_2$ is empty or $t \leq r$, say $t = r$. Set $U=(T_1 \setminus \{e_1\}) \cup \{c_1, \dots, c_{r-1}\}$ and let $T$ be the cut set obtained from $U$ by Lemma \ref{L.CutsetContained}. Clearly $T \subset T_1 \cup T_2$, $e_1 \in T_1 \setminus T_2$, $c_r \in T_2 \setminus T_1$ and $e_1, c_r \notin T$. Moreover, by construction $W$ is a transversal of $G - U$ and by Lemma \ref{L.CutsetContained} (i), there is a transversal of $G - T$ containing $W$; hence, $T$ is the cut set we were looking for.

Thus, from now on we assume that $e_1$ reconnects $B_1, \dots, B_r$ and another component, say $B_{r+1}$. For the same reason, we suppose that $e_j$ reconnects $B_1, \dots, B_r$ and another component, say $B_{r+j}$, for every $j \in [r]$. Note that, except for $B_1, \dots, B_r$, the other components reconnected by $e_i$ and $e_j$, with $i\neq j$, are different because $e_i$ and $e_j$ are in different components of $G-T_2$ by \eqref{Eq.e_jInComponentsOfG-T_2}. It follows that the connected components of $G - T_2$ are the subgraphs induced by $G$ on the vertices $V(B_1 - \{c_1\}), \dots, V(B_r - \{c_r\}), V(B_{r+1}) \cup \{e_1\}, \dots, V(B_{2r}) \cup \{e_r\}$ and $B_{2r+1}, \dots, B_{|W|}$.
In the rest of the proof we will show that this assumption will lead to a contradiction, showing that this case cannot occur.

Notice that $S \neq \emptyset$. In fact, the connected components of $G - (T_1 \cap T_2)$ are $B_{2r +1}, \dots, B_{|W|}$ and another component induced by $G$ on the vertices $V(B_1) \cup \dots \cup V(B_{2r}) \cup (T_1 \setminus T_2)$. Then $c_G(T) = c_G(T_1 \cap T_2) = |W| - 2r +1$ by Lemma \ref{L.CutsetContained} (iii) because $T_1 \cap T_2 \subseteq T_1 \in \CC(G)$. Moreover, $|T|+1 = |T_2| - r - s + 1 = |W| - r - s$ and the unmixedness of $J_G$ implies that $s = r - 1$. Since $r \geq 2$, it follows that $s = |S| > 0$, hence $S \neq \emptyset$.

Let $v \in S = (T_1 \cap T_2) \setminus T$. Since $v \in T_1$, it reconnects at least two components of $G - T_1$ and they are all among $B_1, \dots, B_{2r}$ by Lemma \ref{L.CutsetContained} (ii). We want to prove that in $G - T_1$ the vertex $v$ reconnects exactly one component among $B_1, \dots, B_r$ and one among $B_{r+1}, \dots, B_{2r}$. Assume first that $v$ reconnects at least two components among $B_1, \dots, B_r$, say $B_1,\dots,B_h$. Let $z$ be the first vertex among $c_1,\dots, c_h, v$ that appears in the ordering of the elements of $T_2 \setminus T$, i.e., there exists a cut set $H = T \cup \{g_1, \dots, g_a, z\}$ not containing the other $h$ vertices of $\{c_1,\dots,c_h,v\} \setminus \{z\}$. It is easy to see that $z$ does not reconnect any components in $G - H$ and this yields a contradiction: this can be seen by looking at Figure \ref{F.graphGminusH}, where the dashed edges are edges in $G$ (note that there could be other edges between $v$ and $c_1,\dots,c_h$, and the vertices $\{g_{a+1},\dots,g_{r+s}\} \setminus \{z\}$ are not drawn). Thus, $v$ cannot reconnect more than one component among $B_1, \dots, B_r$. By a similar argument and using strong accessibility on $T \subset T_1$, it also follows that in $G - T_2$ the vertex $v$ cannot reconnect more than one component among the subgraphs induced on each of the sets $V(B_{r+1}) \cup \{e_1\}, \dots, V(B_{2r}) \cup \{e_r\}$. Therefore in $G-T_1$ the vertex $v$ cannot reconnect more than one component among $B_{r+1}, \dots, B_{2r}$ and hence in $G - T_1$ it reconnects exactly one component among $B_1, \dots, B_r$ and one among $B_{r+1}, \dots, B_{2r}$.

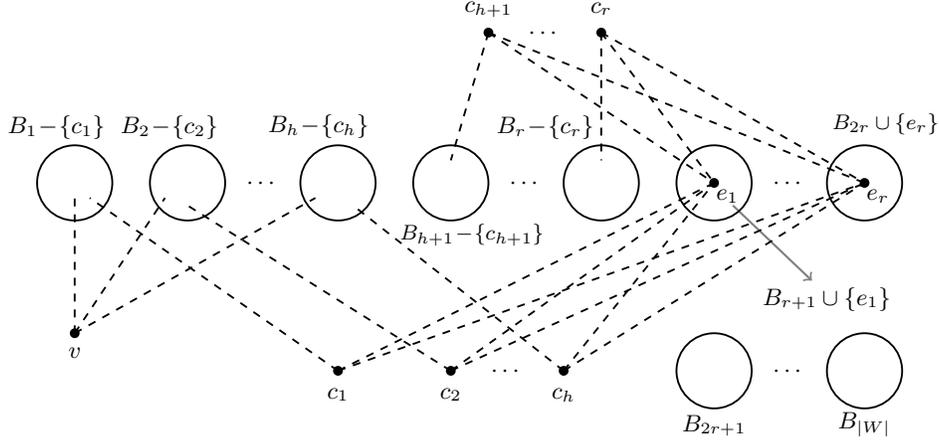
\begin{figure}[ht!]
\begin{tikzpicture}[scale=1]
\node[label={below:{\small $v$}}] (a) at (1,-2) {};
\node[label={below:{\small $c_1$}}] (b) at (4.5,-2.5) {};
\node[label={below:{\small $c_2$}}] (c) at (6,-2.5) {};
\node[label={[label distance=-3mm]below:{\small $\cdots$}}, draw=white, fill=white] (i) at (6.75,-2.5) {};
\node[label={below:{\small $c_h$}}] (d) at (7.5,-2.5) {};
\node[label={[label distance=-1mm]above:{\small $c_{h+1}$}}] (e) at (6.5,2) {};
\node[label={[label distance=-3mm]below:{\small $\cdots$}}, draw=white, fill=white] (i) at (7.25,2) {};
\node[label={above:{\small $c_r$}}] (f) at (8,2) {};
\node[label={[label distance=0mm]below right:{\small $e_1$}}] (g) at (9.5,0) {};
\node[label={[label distance=0mm]below right:{\small $e_r$}}] (h) at (11.5,0) {};

\node[label={[label distance=-4mm]above left:{\small $B_1\! -\! \{c_1\}$}}, draw=white, fill=white] (i) at (1,0.5) {};
\node[label={[label distance=-4mm]above left:{\small $B_2\! -\! \{c_2\}$}}, draw=white, fill=white] (i) at (2.5,0.5) {};
\node[label={[label distance=-3mm]below:{\small $\cdots$}}, draw=white, fill=white] (i) at (3.5,0) {};
\node[label={[label distance=-4mm]above left:{\small $B_h\! -\! \{c_h\}$}}, draw=white, fill=white] (i) at (4.5,0.5) {};

\node[label={[label distance=-5mm]below right:{\footnotesize $B_{h+1}\! -\! \{c_{h+1}\}$}}, draw=white, fill=white] (i) at (5.9,-0.3) {};
\node[label={[label distance=-3mm]below:{\small $\cdots$}}, draw=white, fill=white] (i) at (7,0) {};
\node[label={[label distance=-4mm]above left:{\small $B_r\! -\! \{c_r\}$}}, draw=white, fill=white] (i) at (7.5,0.5) {};

\draw[draw=gray, thick, ->] (9.75,-0.3) -- (10.8,-1.3);

\node[label={[label distance=-4mm]below:{\small $B_{r+1} \cup \{e_1\}$}}, draw=white, fill=white] (i) at (11,-1) {};
\node[label={[label distance=-3mm]below:{\small $\cdots$}}, draw=white, fill=white] (i) at (10.5,0) {};
\node[label={[label distance=-4mm]above right:{\footnotesize $B_{2r} \cup \{e_r\}$}}, draw=white, fill=white] (i) at (11.5,0.5) {};

\node[label={[label distance=-3mm]below:{\small $B_{2r+1}$}}, draw=white, fill=white] (i) at (9.5,-3) {};
\node[label={[label distance=-3mm]below:{\small $\cdots$}}, draw=white, fill=white] (i) at (10.5,-2.5) {};
\node[label={[label distance=-2.5mm]below:{\small $B_{|W|}$}}, draw=white, fill=white] (i) at (11.5,-3) {};

\draw (1,0) circle (0.5cm);
\draw (2.5,0) circle (0.5cm);
\draw (4.5,0) circle (0.5cm);
\draw (6,0) circle (0.5cm);
\draw (8,0) circle (0.5cm);
\draw (9.5,0) circle (0.5cm);
\draw (11.5,0) circle (0.5cm);
\draw (9.5,-2.5) circle (0.5cm);
\draw (11.5,-2.5) circle (0.5cm);

\draw[dashed] (1,-2) -- (1,-0.2);
\draw[dashed] (1,-2) -- (2.2,-0.2);
\draw[dashed] (1,-2) -- (4.2,-0.2);

\draw[dashed] (4.5,-2.5) -- (1.2,-0.2);
\draw[dashed] (4.5,-2.5) -- (9.5,0);
\draw[dashed] (4.5,-2.5) -- (11.5,0);

\draw[dashed] (6,-2.5) -- (2.5,-0.3);
\draw[dashed] (6,-2.5) -- (9.5,0);
\draw[dashed] (6,-2.5) -- (11.5,0);

\draw[dashed] (7.5,-2.5) -- (4.7,-0.3);
\draw[dashed] (7.5,-2.5) -- (9.5,0);
\draw[dashed] (7.5,-2.5) -- (11.5,0);

\draw[dashed] (6.5,2) -- (6,0.3);
\draw[dashed] (6.5,2) -- (9.5,0);
\draw[dashed] (6.5,2) -- (11.5,0);

\draw[dashed] (8,2) -- (8,0.3);
\draw[dashed] (8,2) -- (9.5,0);
\draw[dashed] (8,2) -- (11.5,0);
\end{tikzpicture}
\caption{A partial representation of the graph $G - H$}\label{F.graphGminusH}
\end{figure}

Now we prove that $r=2$. Without loss of generality, up to relabeling, we may assume that in the above ordering $\{g_1, \dots, g_{2r -1}\}$ of the vertices of $S \cup (T_2 \setminus T_1)$, $c_i$ precedes $c_j$ if $i < j$. We claim that in $G - T_1$ every $v \in S$ reconnects $B_r$ to one component among $B_{r+1}, \dots, B_{2r}$. Indeed, suppose that $v$ reconnects $B_i$ with $i < r$ and a component among $B_{r+1}, \dots, B_{2r}$, say $B_{r+1}$. If $c_i$ precedes $v$ in the above order, then there is a cut set $T \cup \{g_1, \dots, g_a, c_i\}$. Note that the graph $G - (T \cup \{g_1, \dots, g_a, c_i\})$ has a connected component containing $V(B_i - \{c_i\}) \cup \{v\}$, $V(B_r) \cup \{c_r\}$, $V(B_{r+1}) \cup \{e_1\}, \dots, V(B_{2r}) \cup \{e_r\}$ and $\{g_{a+1},\dots,g_{2r-1}\} \setminus \{c_i\}$, because the vertices $g_j$ are either in $\{c_{i+1}, \dots, c_r\}$ or in $S$. Since $c_i$ is only adjacent to vertices of $B_i$, to $\{e_1,\dots,e_r\}$ and possibly to some of $g_{a+1},\dots,g_{2r-1}$, it does not reconnect any components in $G - (T \cup \{g_1, \dots, g_a, c_i\})$, a contradiction. On the other hand, if $v$ precedes $c_i$, we conclude in the same way and the claim follows. By symmetry and using strong accessibility on $T \subset T_2$, it also follows that every $v \in S$ reconnects $B_r$ and a fixed component among $B_{r+1}, \dots, B_{2r}$, say $B_{2r}$. Moreover, Remark \ref{R.union} (3) ensures that in a cut set there is at most one vertex that reconnects two fixed components; thus $|S| = s = 1$, and $r-1 = s = 1$ implies $r=2$.

Hence, in $G - T_1$ we have $B_1$ and $B_2$ containing $c_1$ and $c_2$ respectively. Moreover, $e_1$ reconnects $B_1, B_2$, and $B_3$, while $e_2$ reconnects $B_1, B_2$, and $B_4$. Furthermore, $|S|=1$ and by the previous argument $v \in S$ reconnects $B_2$ and $B_4$ and it is the unique vertex of $T_1 \cap T_2$ that reconnects only components among $B_1, B_2, B_3$, and $B_4$. Now consider $H = T_1 \setminus \{e_1\}$. The connected components of $G - H$ are $B_4, \dots, B_{|W|}$, and the component $B$ induced by $G$ on $V(B_1) \cup V(B_2) \cup V(B_3) \cup \{e_1\}$. Then, both $e_2$ and $v$ reconnect $B$ to $B_4$, while all the other vertices in $T_1 \cap T_2$ (if any) reconnect some components of $G - H$ because in $G - T_1$ they reconnect some components that are not all included among $B_1, B_2$, and $B_3$. Then $H$ is a cut set, but $e_1$ reconnects more than two connected components of $G - T_1$; this is a contradiction by Remark \ref{R.union} (2) and it concludes the proof.
\end{proof}

We are now ready to prove Conjecture \ref{Conj.S2_accessible}.

\begin{theorem}\label{T.S2BinomialEdgeIdeals}
Let $G$ be a graph. Then, $R/J_G$ satisfies Serre's condition $(S_2)$ if and only if $G$ is accessible.
\end{theorem}

\begin{proof}
First of all, notice that if $G = \sqcup_{i=1}^r G_i$, where $G_i$ are the connected components of $G$, then $R/J_G$ satisfies Serre's condition $(S_2)$ if and only if the same holds for $R_i/J_{G_i}$ for every $i$, where $R_i = K[x_j, y_j : j \in V(G_i)]$. This follows by the fact that $R/J_G \cong R_1/J_{G_1} \otimes_K \cdots \otimes_K R_r/J_{G_r}$. Thus, it is enough to prove the claim for $G$ connected.

The first implication was proved in \cite[Theorem 2]{LMRR23}. Hence, assume $G$ is accessible. Notice that by \cite[Corollary 2.11]{CV20}, it is enough to prove that $R/\gin(J_G) = R/I_{\Delta_G}$ satisfies Serre’s condition $(S_2)$, where the facets of $\Delta_G$ are described in Lemma \ref{facets}. To do that we are going to apply Theorem \ref{T.characterizationS2}. Suppose that $G$ is an accessible graph on the vertex set $[n]$. If $n \leq 2$, then $R/J_G$ is Cohen-Macaulay and there is nothing to show. Hence, we may assume $n \geq 3$. For every subset $\mbf a \subset [n]$, we set $\mbf{x}_{\mbf a} = \prod_{i \in \mbf a}x_i$ and $\mbf{y}_{\mbf a} = \prod_{i \in \mbf a}y_i$. Let $H = \mbf{y}_{\mbf a} \mbf{x}_{\mbf b}$ be a proper face of $\Delta_G$ for some ${\mbf a}, {\mbf b} \subset [n]$; we only need to show that $\link_{\Delta_G}(H)$ is connected or zero-dimensional.

If $F(S, W) = \mbf{y}_{\mbf c} \mbf{x}_{\mbf d}$ is a facet of $\Delta_G$ containing $H$, we denote by $\widetilde F(S,W)$ the facet of $\link_{\Delta_G}(H)$ obtained from $F(S, W)$, i.e., $\widetilde F(S, W) = \mbf{y}_{\mbf{c} \setminus \mbf{a}} \mbf{x}_{\mbf{d} \setminus \mbf{b}}$.

First we treat some simple cases:
\begin{itemize}
\item If ${\mbf a} = {\mbf b} = \emptyset$, then $H = \emptyset$ and $\link_{\Delta_G}(H) = \Delta_G$, which is connected by Proposition \ref{P.gin_cone}.
\item If $\mbf{a} = \emptyset$, then $H = \mbf{x}_{\mbf b}$ for some $\mbf{b} \subset [n]$ and hence the facets of $\link_{\Delta_G}(H)$ have the form $\mbf{y}_\mbf{c} \mbf{x}_\mbf{d} $, with $\mbf{b} \cup \mbf{d} \subset \mbf{c}$. It follows that $\link_{\Delta_G}(H)$ is connected since all facets of $\link_{\Delta_G}(H)$ contain $\mbf{y}_{\mbf b}$.
\item If $\mbf{b} = \emptyset$, then $H = \mbf{y}_\mbf{a}$ for some $\mbf{a} \subset [n]$. Notice that, if $|\mbf{a}| = n$, then
    \[
    \link_{\Delta_G}(H) = \{\widetilde F(\emptyset, \{j\}) : j = 1, \dots, n\} = \{x_1, \dots, x_n\}
    \]
    is zero-dimensional. Hence, we may assume $|\mbf{a}| < n$; in this case, $\widetilde F(\emptyset, \{j\}) = \mbf{y}_{[n] \setminus \mbf{a}}x_j$ is a facet of $\link_{\Delta_G}(H)$ for every $j \in [n]$. These facets are all connected through $\mbf{y}_{[n] \setminus \mbf{a}}$. Since any other facet of $\link_{\Delta_G}(H)$ contains some $x_j$, $\link_{\Delta_G}(H)$ is connected.
\end{itemize}

Let $H = \mbf{y}_\mbf{a} \mbf{x}_\mbf{b}$, with $|\mbf{a}| \geq 1$, $|\mbf{b}| \geq 1$ and $|\mbf{a}| + |\mbf{b}| < n$, otherwise $\dim(\link_{\Delta_G}(H)) = 0$ because $\dim(\Delta_G) = n$. Consider two facets of $\Delta_G$ containing $H$:
\[
F_1 = F(T_1, W_1) = (\mbf{y}_\mbf{a} \mbf{x}_\mbf{b}) \mbf{y}_\mbf{c} \mbf{x}_\mbf{d} \text{\ \ and\ \ } F_2 = F(T_2, W_2) = (\mbf{y}_\mbf{a} \mbf{x}_\mbf{b}) \mbf{y}_\mbf{e} \mbf{x}_\mbf{f},
\]
where $T_1, T_2 \in \CC(G)$ and $W_1, W_2$ are transversals of $G - T_1$ and of $G - T_2$, respectively. In particular, $\widetilde F_1 = \mbf{y}_\mbf{c} \mbf{x}_\mbf{d}$ and $\widetilde F_2 = \mbf{y}_\mbf{e} \mbf{x}_{\mbf f}$ are facets of $\link_{\Delta_G}(H)$, with $\mbf{c}, \mbf{d}, \mbf{e}, \mbf{f} \subset [n]$.

\smallskip
\textbf{Claim}: The facets $\widetilde F_1$ and $\widetilde F_2$ are connected through other facets of $\link_{\Delta_G}(H)$, i.e., $\widetilde F_1 = \widetilde L_0, \widetilde L_1, \dots, \widetilde L_s = \widetilde F_2$, where $\widetilde L_i$ are facets of $\link_{\Delta_G}(H)$ and $\widetilde L_i \cap \widetilde L_{i-1} \neq \emptyset$ for every $i=1,\dots,s$.

\smallskip
First of all, we may assume $\mbf{c} \cap \mbf{e} = \mbf{d} \cap \mbf{f} = \emptyset$, otherwise $\widetilde F_1$ and $\widetilde F_2$ are connected. Moreover, $|\mbf{c}| + |\mbf{d}| = |\mbf{e}| + |\mbf{f}|$ since $\Delta_G$ is pure by Lemma \ref{facets} and we assume $|\mbf{c}| + |\mbf{d}|>1$ otherwise $\link_{\Delta_G}(H)$ has dimension $0$.

Notice that $W_1 = \mbf{b} \cup \mbf{d}$, $W_2 = \mbf{b} \cup \mbf{f}$, $\mbf{a} = [n] \setminus (T_1 \cup T_2)$, $\mbf{c} = T_2 \setminus T_1$ and $\mbf{e} = T_1 \setminus T_2$. By the unmixedness of $J_G$, we have that $|W_1| = |\mbf{b}| + |\mbf{d}| = |T_1| + 1$ and $|W_2| = |\mbf{b}| + |\mbf{f}| = |T_2| + 1$.

\textbf{Case 1)}: Assume $\mbf{e} = \emptyset$ (or similarly $\mbf{c} = \emptyset$).

In this case, $|\mbf{f}| \geq 2$ and $T_1 \subset T_2$. If also $\mbf{c} = \emptyset$, then $T_1 = T_2$ and we denote this cut set by $T$. Therefore,\break $|\mbf{d}| = |\mbf{f}| \geq 2$. Consider $d \in \mbf{d}$ and $f \in \mbf{f}$ belonging to the same connected component of $G - T$. Then\break $W = \mbf{b} \cup (\mbf{d} \setminus \{d\}) \cup \{f\}$ is a transversal of $G - T$ and $L = F(T, W)$ is a facet of $\Delta_G$ such that $H \subset L$, hence $\widetilde L = \mbf{x}_{(\mbf{d} \setminus \{d\}) \cup \{f\}}$ is a facet of $\link_{\Delta_G}(H)$. Since $|\mbf{d}| \geq 2$, $\widetilde F_1$ is connected to $\widetilde{L}$ through $\mbf{x}_{\mbf{d} \setminus \{d\}}$, whereas $\widetilde{L}$ is connected to $\widetilde F_2$ through $x_f$:
\[
\widetilde F_1 = \mbf{x}_\mbf{d}\ \frac{\ \ \mbf{x}_{(\mbf{d} \setminus \{d\})}\ \ }{ }\ \widetilde L = \mbf{x}_{(\mbf{d} \setminus \{d\}) \cup \{f\}}\ \frac{\ \ x_f\ \ }{ }\ \widetilde F_2 = \mbf{x}_\mbf{f}.
\]

Assume now $\mbf{c} \neq \emptyset$. Since $T_1 \subsetneq T_2$, $\CC(G)$ is strongly accessible by Theorem \ref{T.stronglyAccessible} and hence there exists\break $c \in \mbf{c}=T_2 \setminus T_1$ such that $T_2 \setminus \{c\} \in \CC(G)$. In particular, $c$ reconnects exactly two components of $G - T_2$ by Remark \ref{R.union} (2) and it cannot reconnect two components containing two elements of $\mbf{b}$, because $T_1 \subset T_2 \setminus \{c\}$ and in $G - T_1$ the vertices of $\mbf{b}$ belong to pairwise distinct components. It follows that $c$ reconnects at least one component containing some $f \in \mbf{f}$. Then $W = \mbf{b} \cup (\mbf{f} \setminus \{f\})$ is a transversal of $G - (T_2 \setminus \{c\})$ and $L = F(T_2 \setminus \{c\}, W)$ is a facet of $\Delta_G$ with $H \subset L$, hence $\widetilde L = y_c \mbf{x}_{\mbf{f} \setminus \{f\}}$ is a facet of $\link_{\Delta_G}(H)$. Notice that $\widetilde F_1$ is connected to $\widetilde{L}$ through $y_c$ and $\widetilde{L}$ is connected to $\widetilde F_2$ through $\mbf{x}_{\mbf{f} \setminus \{f\}}$ (recall that $|\mbf{f}| \geq 2$):
\[
\widetilde F_1 = \mbf{y}_\mbf{c} \mbf{x}_\mbf{d}\ \frac{\ \ y_c\ \ }{ }\ \widetilde L = y_c \mbf{x}_{\mbf{f} \setminus \{f\}}\ \frac{\ \ \mbf{x}_{\mbf{f} \setminus \{f\}}\ \ }{ }\ \widetilde F_2 = \mbf{x}_\mbf{f}.
\]

\textbf{Case 2)}: Let $\mbf{c} \neq \emptyset$ and $\mbf{e} \neq \emptyset$.

\smallskip
\textbf{Case 2.1)}: Assume $\mbf{d} \neq \emptyset$ and $\mbf{f} \neq \emptyset$.

Let $d \in \mbf{d}$ and $f \in \mbf{f}$. Consider a cut set $T \subset T_1 \cup T_2$ obtained by using Lemma \ref{L.CutsetContained}. We claim that we may assume $d \notin T$. Indeed, assume $d \in T$ and let $B$ be the connected component of $G - T_1$ containing $d$. If $V(B) \subset T \subset T_1 \cup T_2$, then $d$ does not reconnect any connected components in $G - (T_1 \cup T_2)$ and we may assume that $d \notin T$ by Lemma \ref{L.CutsetContained} (iv). On the other hand, if $V(B) \not \subset T$, let $v \in V(B) \setminus T$. Notice that $W = (W_1 \setminus \{d\}) \cup \{v\}$ is a transversal of $G - T_1$, hence $L = F(T_1, W)$ is a facet of $\Delta_G$ with $H \subset L$. Since $\widetilde{F}_1$ and $\widetilde{L}$ are connected through $\mbf{y}_\mbf{c}$, we may replace $d$ with $v$ and $v \notin T$. Similarly, we may assume $f \notin T$. The vertices of $\mbf{b} \cup \{d\}$  are in different connected components of $G - T$, because they are in different components of $G - T_1$. The same holds for the vertices of $\mbf{b} \cup \{f\}$. Let $C$ be the connected component of $G - T$ containing $d$.

If $f \notin C$, then there exists $\mbf{w} \subset [n]$ such that $Z = \mbf{b} \cup \mbf{w} \cup \{d,f\}$ is a transversal of $G - T$. In this case, $L = F(T, Z)$ is a facet of $\Delta_G$ such that $H \subset L$, because $T \subset T_1 \cup T_2$. Thus, $\widetilde{F}_1$ is connected to $\widetilde{L}$ through $x_d$ and $\widetilde{L}$ is connected to $\widetilde{F}_2$ through $x_f$:
\[
\widetilde F_1 = \mbf{y}_\mbf{c} \mbf{x}_\mbf{d}\ \frac{\ \ x_d\ \ }{ }\ \widetilde L = \mbf{y}_{[n] \setminus (\mbf{a} \cup T)} \mbf{x}_{\mbf{w} \cup \{d,f\}}\ \frac{\ \ x_f\ \ }{ }\ \widetilde F_2 = \mbf{y}_\mbf{e} \mbf{x}_\mbf{f}.
\]

If $f \in C$, then there exists $\mbf{w} \subset [n]$ such that $Z_1 = \mbf{b} \cup \mbf{w} \cup \{d\}$ and $Z_2 = \mbf{b} \cup \mbf{w} \cup \{f\}$ are transversals of $G - T$. Notice that if $T = T_1 \cup T_2$, then $\mbf{w} \neq \emptyset$. In fact, $|\mbf{w} \cup \{d\}| = |Z_1 \setminus \mbf{b}| = |Z_1| - |\mbf{b}| = |T_1 \cup T_2| + 1 - |\mbf{b}| = n - |\mbf{a}| + 1 - |\mbf{b}| > 1$. As before, $L_1 = F(T, Z_1)$ and $L_2 = F(T, Z_2)$ are facets of $\Delta_G$, both containing $H$. Hence, $\widetilde{F}_1$ is connected to $\widetilde{L}_1$ through $x_d$, $\widetilde{L}_2$ is connected to $\widetilde{F}_2$ through $x_f$, and $\widetilde{L}_1$ is connected to $\widetilde{L}_2$ through $\mbf{y}_{[n] \setminus (\mbf{a} \cup T)}$ (if $T \neq T_1 \cup T_2$) or through $\mbf{x}_\mbf{w}$ (if $T = T_1 \cup T_2$):
\[
\widetilde F_1 = \mbf{y}_\mbf{c} \mbf{x}_\mbf{d}\ \frac{\ \ x_d\ \ }{ }\ \widetilde L_1 = \mbf{y}_{[n] \setminus (\mbf{a} \cup T)} \mbf{x}_{\mbf{w} \cup \{d\}}\ \frac{\ \ \mbf{y}_{[n] \setminus (\mbf{a} \cup T)} \text{ or } \mbf{x}_\mbf{w}\ \ }{ }\ \widetilde L_2 = \mbf{y}_{[n] \setminus (\mbf{a} \cup T)} \mbf{x}_{\mbf{w} \cup \{f\}}\ \frac{\ \ x_f\ \ }{ }\ \widetilde F_2 = \mbf{y}_\mbf{e} \mbf{x}_\mbf{f}.
\]

\textbf{Case 2.2)}: Assume $\mbf{d} = \emptyset$ and $\mbf{f} \neq \emptyset$ (or similarly $\mbf{d} \neq \emptyset$ and $\mbf{f}=\emptyset$).

In this case, we have $\mbf{b} = W_1 \subsetneq W_2 = \mbf{b} \cup \mbf{f}$ and $|\mbf{c}| = |\mbf{e}| + |\mbf{f}| \geq 2$. Let $f \in \mbf{f} = W_2 \setminus W_1$. If $T_1 \cup \{c\} \in \CC(G)$ for some $c \in \mbf{c}$, then there exists $j \notin \mbf{b}$ such that $\mbf{b} \cup \{j\}$ is a transversal of $G - (T_1 \cup \{c\})$. Then $L = F(T_1 \cup \{c\}, \mbf{b} \cup \{j\})$ is a facet of $\Delta_G$ such that $H \subset L$. Clearly, $\widetilde F_1$ is connected to $\widetilde{L}$ through $\mbf{y}_{\mbf{c} \setminus \{c\}}$ and $\widetilde{L}$ is connected to $\widetilde F_2$ by \textbf{Case 2.1)}:
\[
\widetilde F_1 = \mbf{y}_\mbf{c}\ \frac{\ \ \mbf{y}_{\mbf{c} \setminus \{c\}}\ \ }{ }\ \widetilde L = \mbf{y}_{\mbf{c} \setminus \{c\}} x_j\ \frac{\ \ \textbf{Case 2.1}\ \ }{ }\ \widetilde F_2 = \mbf{y}_\mbf{e} \mbf{x}_\mbf{f}.
\]
Therefore, we may assume that $T_1 \cup \{c\} \notin \CC(G)$ for every $c \in \mbf{c} = T_2 \setminus T_1$.

If $f \notin T_1$, by Lemma \ref{L.firstLemmaMainTheorem}, there exist a cut set $T \subset T_1 \cup T_2$ of $G$, with $T_1 \not \subset T$, $T_2 \setminus T_1 \not \subset T$, and $\mbf{b} \cup \mbf{w}_1$ is a transversal of $G - T$ for some $\mbf{w}_1 \subset [n]$. Moreover, there exists $c \in T_2 \setminus T_1$, $c \notin T$, such that $T \cup \{c\} \in \CC(G)$ and $\mbf{b} \cup \{f\} \cup \mbf{w}_2$ is a transversal of $G - (T \cup \{c\})$ for some $\mbf{w}_2 \subset [n]$. It follows that $L_1 = F(T, \mbf{b} \cup \mbf{w}_1)$ and $L_2 = F(T \cup \{c\}, \mbf{b} \cup \{f\} \cup \mbf{w}_2)$ are facets of $\Delta_G$, both containing $H$. Since $T_1 \not \subset T$, there exists $s \in T_1 \setminus T$, hence $s \notin \mbf{a}$. Thus $\widetilde{F}_1$ is connected to $\widetilde{L}_1$ through $y_c$, $\widetilde{L}_1$ is connected to $\widetilde{L}_2$ through $y_s$ and $\widetilde{L}_2$ is connected to $\widetilde{F}_2$ through $x_f$:
\[
\widetilde F_1 = \mbf{y}_\mbf{c}\ \frac{\ \ y_c\ \ }{ }\ \widetilde L_1 = \mbf{y}_{[n] \setminus (\mbf{a} \cup T)} \mbf{x}_{\mbf{w}_1}\ \frac{\ \ y_s\ \ }{ }\ \widetilde L_2 = \mbf{y}_{[n] \setminus (\mbf{a} \cup T \cup \{c\})} \mbf{x}_{\mbf{w}_2 \cup \{f\}}\ \frac{\ \ x_f\ \ }{ }\ \widetilde F_2 = \mbf{y}_\mbf{e} \mbf{x}_\mbf{f}.
\]

Assume $f \in T_1$ and suppose first that $N_G(f) \not \subset T_1 \cup T_2$. Let $g \notin T_1 \cup T_2$ such that $g \in N_G(f)$. Then $W = (W_2 \setminus \{f\}) \cup \{g\}$ is a transversal of $G - T_2$ and $L = F(T_2, W)$ is a facet of $\Delta_G$ with $H \subset L$. Thus $\widetilde{F}_2$ is connected to $\widetilde{L}$ through $\mbf{y}_\mbf{e}$ and $\widetilde{L}$ to $\widetilde{F}_1$ as above, because $g \notin T_1$. Suppose now that $N_G(f) \subset T_1 \cup T_2$. By Lemma \ref{L.secondLemmaMainTheorem}, there exists a cut set $T \subset T_1 \cup T_2$ of $G$, with $T_1 \setminus T_2 \not \subset T$, $T_2 \setminus T_1 \not \subset T$ and $\mbf{w} \subset [n]$ such that $\mbf{b} \cup \mbf{w}$ is a transversal of $G - T$. Hence $L = F(T, \mbf{b} \cup \mbf{w})$ is a facet of $\Delta_G$ containing $H$. Let $c \in T_2 \setminus T_1$ with $c \notin T$ and let $e \in T_1 \setminus T_2$ with $e \notin T$. Thus, $\widetilde{F}_1$ is connected to $\widetilde{L}$ through $y_c$ and $\widetilde{L}$ is connected to $\widetilde{F}_2$ through $y_e$:
\[
\widetilde F_1 = \mbf{y}_\mbf{c}\ \frac{\ \ y_c\ \ }{ }\ \widetilde L = \mbf{y}_{[n] \setminus (\mbf{a} \cup T)} \mbf{x}_{\mbf{w}}\ \frac{\ \ y_e\ \ }{ }\ \widetilde F_2 = \mbf{y}_\mbf{e} \mbf{x}_\mbf{f}.
\]

\textbf{Case 2.3)} Assume $\mbf{d} = \mbf{f} = \emptyset$.

We have $W_1 = W_2 = \mbf{b}$ and $|\mbf{c}| = |\mbf{e}| \geq 2$, hence $|T_1| = |T_2|$. We may assume that there is no cut set $S \subset T_1 \cup T_2$ and with a transversal $W$ for which $\mbf{b} \subsetneq W$. In fact, if such a cut set exists, then $L = F(S,W)$ would be a facet of $\Delta_G$ containing $H$; in particular, both $\widetilde{F}_1$ and $\widetilde{F}_2$ would be connected to $\widetilde{L}$ by \textbf{Case 1)} or by \textbf{Case 2.2)}.

By Lemma \ref{L.thirdLemmaMainTheorem}, there exists a cut set $T \subset T_1 \cup T_2$ of $G$ such that $T_1 \setminus T_2 \not \subset T$, $T_2 \setminus T_1 \not \subset T$ and $\mbf{b} \cup \mbf{w}$ is a transversal of $G - T$ for some $\mbf{w} \subset [n]$. Hence $L = F(T, \mbf{b} \cup \mbf{w})$ is a facet of $\Delta_G$ containing $H$. Let $c \in T_2 \setminus T_1$ such that $c \notin T$ and $e \in T_1 \setminus T_2$ such that $e \notin T$. Hence, $\widetilde{F}_1$ is connected to $\widetilde{L}$ through $y_c$ and $\widetilde{L}$ is connected to $\widetilde{F}_2$ through $y_e$:
\[
\widetilde F_1 = \mbf{y}_\mbf{c}\ \frac{\ \ y_c\ \ }{ }\ \widetilde L = \mbf{y}_{[n] \setminus (\mbf{a} \cup T)} \mbf{x}_{\mbf{w}}\ \frac{\ \ y_e\ \ }{ }\ \widetilde F_2 = \mbf{y}_\mbf{e}.
\]

We conclude that all the non-zero dimensional links of $\Delta_G$ are connected, thus $R/I_{\Delta_G} = R/\gin(J_G)$ satisfies Serre's condition $(S_2)$ by Theorem \ref{T.characterizationS2}. This is sufficient to conclude that $R/J_G$ satisfies Serre's condition $(S_2)$ by \cite[Corollary 2.11]{CV20}.
\end{proof}

\section*{Acknowledgements}
The second author was supported by the Deutsche Forschungsgemeinschaft (DFG, German Research Foundation) – project number 454595616.
Part of this paper was written while the first and the fourth authors visited the Discrete Mathematics and Topological Combinatorics research group at Freie Universit\"{a}t Berlin. They want to express their thanks for the hospitality.

\end{document}